\documentclass[a4paper,11pt]{article}
\usepackage{eurosym}
\usepackage{framed,amssymb,amsmath,amsthm,bm,color,soul,url}
\usepackage[normalem]{ulem}
\usepackage{tikz}
\usepackage{authblk}

\usetikzlibrary{fit,chains}
\newtheorem{theorem}{Theorem}[section]
\newtheorem{proposition}[theorem]{Proposition}
\newtheorem{remark}[theorem]{Remark}
\newtheorem{lemma}[theorem]{Lemma}
\newtheorem{corollary}[theorem]{Corollary}
\newtheorem{definition}[theorem]{Definition}

\usepackage{graphicx}

\title{A family of systems including the Herschel-Bulkley fluid equations}
\author[ ]{Nikolai V. Chemetov\thanks{Department of Computing and Mathematics, University of S{\~a}o Paulo. 14040-901 Ribeir{\~a}o Preto - SP, Brazil. Email:
nvchemetov@gmail.com}}
\author[ ]{Marcelo M. Santos\thanks{IMECC-Universidade Estadual de Campinas (UNICAMP). Rua S\'ergio Buarque de Holanda, 651. 13083-859 Campinas, SP, Brazil. Email: mmsantos@unicamp.br}} 
\affil{ }

\date{ }

\begin{document}

\maketitle

\begin{abstract}
We analyze the Navier-Stokes equations for incompressible fluids with the {\lq\lq}viscous stress tensor{\rq\rq} $\mathbb{S}$ in a family which includes the Bingham model for viscoplastic fluids (more generally, the Herschel-Bulkley model).  $\mathbb{S}$ is  the subgradient of a convex  potential $V=V(x,t,X)$, allowing  that $V$ can depend on the space-time variables $(x,t)$. The potential has its one-sided directional derivatives $V'(X,X)$ uniformly bounded from below and above by a $p$-power function of the matrices $X$. For $p\geqslant 2.2$ we solve an initial boundary value problem for those fluid systems, in a bounded region in $\mathbb{R}^3$. We take a  nonlinear boundary condition, which encompasses the Navier friction/slip boundary condition.
\end{abstract}


\smallskip 


\section{Introduction}

\label{Introduction}

The Navier-Stokes system
\begin{equation}
\mathbf{v}_{t}+\left(\mathbf{v\cdot\nabla}\right) \mathbf{v}=\mbox{div}\ \mathbb{T},\quad \mbox{div}\,\mathbf{v}=0 
\label{NS}
\end{equation}
models an incompressible fluid with particle velocity $\mathbf{v}=\mathbf{v}(x,t)$ at a position $x$ in a time $t$. The first equation (in fact, a system of equations) expresses the balance of momentum on the particle under internal forces\footnote{For simplicity of presentation, we disregard here external forces.}  represented
by the term $\mbox{div}\ \mathbb{T}$, where $\mathbb{T}$ is the so-called {\em Cauchy stress tensor} $\mathbb{T}$. The internal forces are given by the pressure force, viscous-related forces and forces due to possible rotation of the particle fluids\footnote{For the case of fluids containing rotating particles, we refer to \cite{E,Lukas,ShelRuz2013,T}.}. It is well established that
$$
\mathbb{T}=-{\rm p}\,\mathrm{I}+\mathbb{S},
$$
where $\mathbb{S}$ is the {\lq\lq}viscous part{\rq\rq} of the Cauchy stress tensor (see e.g. \cite{Malek-Necas-Rokyta-Ruzicka}). We shall call $\mathbb{S}$ the \emph{viscous stress tensor} and throughout the paper, to shorten the notation, we set $\mathrm{B}=\nabla \mathbf{v}$.

Here, we 
we are interested in possible effects of rotating particle in the tensor $\mathbb{S}$ as well. Let us remark that in such case the tensor $\mathbb{S}$ depends on both the symmetric and antisymmetric parts of $\nabla\mathbf{v}$ (in contrast to the case where rotations are not taken into account - in this case, as well known, the tensor $\mathbb{S}$ depends only on the symmetric part of $\nabla\mathbf{v}$). See e.g. \cite{E,Lukas,ShelRuz2013,T}. 

As there is a great variety of fluids, the specification of $\mathbb{S}$ is a major issue in the mathematical modeling of fluids. 
Viscoplastic fluids are materials that exhibit a combination of solid-like (plastic) behavior and fluid-like (viscous) behavior under different conditions. They are a type of non-Newtonian fluid that have a yield stress, which means they behave like solids at low stresses but flow like liquids once a certain threshold stress (yield stress) is exceeded. More precisely, at each time the region occupied by the fluid is divided in two parts: one where the stress surpasses the yield stress (resulting in liquid-like flow), and the other with lower stress (where the material behaves like a solid). Some examples of viscoplastic fluids are 
toothpaste, slurries (mixtures of solid particles suspended in a liquid), drilling fluids (see e.g. \cite[Chapter 1]{BR}) and animal blood.

The parts where the fluid exhibits solid-like and liquid-like behaviors at any given time are unknown in fluid models. They can be defined in terms of the yield stress, which is a non negative number, depending (only) on the kind of viscoplastic fluid considered. We shall denoted it by $\tau_{\ast}$. We also shall call the part where the fluid is solid-like {\em plug region}, denote it by $\Omega^*_t$ at a time $t$ and define it at as being the set
$$
\Omega^*_t:=\{x\, ;\, |\mathbb{S}(x,t)|\le \tau_{\ast}\}.
$$

In the simplest (and pionering) model for viscoplastic fluids, the tensor $\mathbb{S}$
is given by
\begin{equation}  
\left\{\begin{array}{rl}
\mathbb{S} & = \mu _{1}\mathrm{B}_{s}+\tau_{\ast }\frac{\mathrm{B}_{s}}{| \mathrm{B}_{s}|}, \quad \mbox{ if }\mathrm{B}
_{s}\neq 0 \\ 
\ |\mathbb{S}|\! & \le \ \tau_{\ast}, \quad \mbox{ if }\mathrm{B}_{s}=0.%
\end{array}\right.
\label{B}
\end{equation}
These relations are a constitutive law for the (viscoplastic) fluid model
known as \emph{Bingham fluid} or {\em Bingham model}. Here, and throughout the paper, for any square real matrix $\mathrm{X}$, we use the decomposition \eqref{decomposition}.\footnote{Throughout this paper, $d\in\mathbb{N}$, $\mathbb{R}^{d\times d}$ denote the space of $d\times d$ real matrices and we use the Euclidean norm
$|\mathrm{X}|:=\sqrt{\mathrm{X}:\mathrm{X}}$, where, for $\mathrm{X}=(X_{ij}),\,\mathrm{Y}=(Y_{ij})\in \mathbb{R}^{d\times d}$, $\mathrm{X}:\mathrm{Y}$ is the standard inner product $\mathrm{X}:\mathrm{Y}:=\sum_{i,j=1}^d X_{ij}Y_{ij}$. We recall the elementary fact that the space of the matrices endowed with this inner product is the direct sum of the spaces of symmetric matrices and antisymmetric matrices, i.e., for any matrices $\mathrm{X},\mathrm{Y}\in \mathbb{R}^{d\times d}$ we have the decomposition 
\begin{equation}
\mathrm{X}=\mathrm{X}_s+\mathrm{X}_a,
\label{decomposition}
\end{equation}
and $\mathrm{X}_s:\mathrm{Y}_a=\mathrm{0}$, where $\mathrm{X}_s$ and $\mathrm{X}_a$ are the symmetric and antisymmetric parts of $\mathrm{X}$, i.e., $\mathrm{X}_s:=\frac{1}{2}(X+X^{trans})$ and $\mathrm{X}_a:=\frac{1}{2}(X-X^{trans})=\mathrm{X}-\mathrm{X}_s$, being $X^{trans}$ the transpose of $\mathrm{X}$. We shall often use the convenient properties: $X_{s}:Y=X_{s}:Y_{s}$, $X_{a}:Y=X_{a}:Y_{a}$, $X_{s}:X=|X_{s}|^{2}$,  $X_{a}:X=|X_{a}|^{2}$,  $|X|^{2}=|X_{s}|^{2}+|X_{a}|^{2}$.}

The coefficient $\mu _{1}$ is a positive constant (viscosity coefficient).\footnote{ We use the subscript ${}_1$ in $\mu _{1}$ to distinguish it from the viscosity coefficient $\mu _{2 }$ we use below.}

Let us observe that the first relation in \eqref{B} can be written as 
$$
\mathbb{S} = (\mu _{1}|\mathrm{B}_{s}|+\tau_{\ast })\frac{\mathrm{B}_{s}}{|\mathrm{B}_{s}|},
$$
then it implies that $|\mathbb{S}|=\mu _{1}|\mathrm{B}_{s}|+\tau_{\ast }>\tau_{\ast }$. Thus, the plug region $\Omega^\ast_t$ for the Bingham model (the model \eqref{B}) is given by  
$$
\Omega^*_t:=\{x\, ;\, \mathrm{B}_{s}=0\}.
$$

Another observation with respect to the relations \eqref{B}, which is a key point in this paper, is that a tensor $\mathbb{S}$ satisfying \eqref{B} has a {\em potential}, i.e., the relations \eqref{B} can be achieved by the single condition
\begin{equation}
\mathbb{S}(x,t)\in \partial V(\mathrm{B}(x,t)),
\label{BV}
\end{equation} 
for each point $(x,t)$, where $V$ is the convex function $V:\mathbb{R}^{d\times d}\to\mathbb{R}$ given by 
$$
V(\mathrm{X})=\frac{\mu _{1}}{2}|\mathrm{X}_{s}|^{2}+\tau _{\ast }|\mathrm{X}_{s}|
$$
and $\partial V(\mathrm{X})$ denotes the subdifferential of $V$ at the point (matrix) $\mathrm{X}$. Indeed, if $\mathrm{X}_s\not=0$ then $V$ is differentiable at $\mathrm{X}$ and $\partial V(\mathrm{X})=\{\nabla V(X)=\mu _{1}\mathrm{X}_{s}+\tau_{\ast }\frac{\mathrm{X}_{s}}{|\mathrm{X}_{s}|}\}$, thus if $\mathrm{B}_s\not=0$ the condition \eqref{BV} is equivalent to the first relation in \eqref{B}. On the other hand, the subdifferential of $V$ at a matrix $\mathrm{X}$ with $\mathrm{X}_s=\mathrm{0}$  is given by $\tau _{\ast }$ times the subdifferential at $\mathrm{X}$ of the function $f(\mathrm{X})=|\mathrm{X}_{s}|$, since the function $\mathrm{X}\mapsto |\mathrm{X}_{s}|^2$ has null derivative at the origin.   By the definition of subdifferential for convex functions\footnote{For the basics and other results on convex analysis we refer to the book \cite{R}.} and $f(\mathrm{0})=0$, we have that a matrix $\mathrm{X}^*$ belongs to the subdifferential of
$f(\mathrm{X})=|\mathrm{X}_{s}|$ at $\mathrm{X}$  with $\mathrm{X}_s=\mathrm{0}$ if and only if
\begin{equation}
|\mathrm{Y}_{s}|\ge \mathrm{X}^*:(\mathrm{Y}-\mathrm{X})
\label{sdXs}
\end{equation}
for all $\mathrm{Y}\in\mathbb{R}^{d\times d}$. Taking $\mathrm{Y}-\mathrm{X}=\mathrm{X}^\ast_s,\, \mathrm{X}^\ast_a$ in \eqref{sdXs} and using the orthogonality $\mathrm{X}_s:\mathrm{Y}_a=0$, for any matrices $\mathrm{X},\mathrm{Y}\in\mathbb{R}^{d\times d}$, we obtain that $|\mathrm{X}^\ast_s|\le 1$ (here we use that $\mathrm{X}_s=\mathrm{0}$)  and $\mathrm{X}^\ast_a=0$. Reciprocally, if $|\mathrm{X}^\ast_s|\le 1$ and $\mathrm{X}^\ast_a=0$ then by Cauchy-Schwarz inequality we have that $\mathrm{X}^\ast:(\mathrm{Y}-\mathrm{X})=\mathrm{X}^\ast_s:\mathrm{Y}_s\le |\mathrm{Y}_s|=f(\mathrm{Y})$ for all $\mathrm{Y}\in\mathbb{R}^{d\times d}$, thus $\mathrm{X}^\ast\in\partial f(\mathrm{X})$. That is, the subdifferential $\partial f(\mathrm{X})$ of the function $f(\mathrm{X})=|\mathrm{X}_{s}|$ at a matrix $\mathrm{X}$ with $\mathrm{X}_s=\mathrm{0}$ is the set $\{\mathrm{X}\in \mathbb{R}^{d\times d}:\ \mathrm{X}_a=0 \mbox{ and }
|\mathrm{X}_{s}|\leqslant 1\}$ (the projection of the unit ball in $\mathbb{R}^{d\times d}$ on the space of symmetric matrices). Consequently, 
\begin{equation}
\partial V(\mathrm{X})=\{\mathrm{X}^*\in \mathbb{R}^{d\times d}:\ \mathrm{X}^*_a=0 \mbox{ and }
|\mathrm{X}^*_{s}|\leqslant {\tau }_{\ast }\},
\label{sdXs0}
\end{equation}
at matrices $\mathrm{X}$ with $\mathrm{X}_s=\mathrm{0}$. Thus, if $\mathrm{B}_s=\mathrm{0}$ and $\mathbb{S}\in\partial V(\mathrm{B})$ then $|\mathbb{S}|\le\tau_{\ast}$.
  
A generalization to the Bingham model is the power law model known as Herschel-Bulkley model (cf. e.g. \cite[(3.1)]{ShelRuz2013})
\begin{equation}
\left\{\begin{array}{rl}
\mathbb{S} & = \ \mu _{1}|\mathrm{B}_{s}|^{p-2}\mathrm{B}_{s}+\tau_{\ast }%
\frac{\mathrm{B}_{s}}{| \mathrm{B}_{s}|}, \quad \mbox{ if }\mathrm{B}%
_{s}\neq 0 \\ 
\ |\mathbb{S}| & \le \ \tau_{\ast}, \quad \mbox{ if }\mathrm{B}_{s}=0,
\end{array}\right.
\label{HB}
\end{equation}
where $p$ is a real number (the power law index). As in the Bingham model (the case $p=2$), the plug zone for \eqref{HB} is the set $\Omega^*_t:=\{x\, ;\, \mathrm{B}_{s}=0\}$ \ (notice that $\mathbb{S}=(\mu _{1}|\mathrm{B}_{s}|^{p-1}+\tau_{\ast })\frac{\mathrm{B}_{s}}{|\mathrm{B}_{s}|}$, if $\mathrm{B}_{s}\not=\mathrm{0}$) and, considering $p\ge1$, it also has a convex potential, which is the function
$$
V(\mathrm{X})=\frac{\mu _{1}}{p}|\mathrm{X}_{s}|^{p}+\tau _{\ast }|\mathrm{X}_{s}|.
$$
We still have here that the subdifferential $\partial V(\mathrm{X})$, at matrices $\mathrm{X}$ with $\mathrm{X}_s=\mathrm{0}$, is given by \eqref{sdXs0}.

In the case of fluid with rotating particles (micropolar viscoplastic fluid), the viscous stress tensor $\mathbb{S}$ depends also on the antisymmetric part of $\mathrm{B}=\nabla\mathbf{v}$, as we already mentioned above. The {\em Cosserat-Bingham model} proposed by Shelukhin and R\r{u}\v{z}i\v{c}ka \cite{ShelRuz2013} has the viscous stress tensor $\mathbb{S}$, in the part of the fluid domain where the stress surpasses the yield stress (the compliment of the plug region), given as a function of $\mathrm{B}_s$ (the symmetric part of $\mathrm{B}$) and 
\begin{equation}
\mathrm{R}=\mathrm{B}_{a}-{\Omega},
\label{R}
\end{equation}
where ${\Omega}={\Omega}(x,t)$ is the \emph{micro-rotational velocity tensor}, which is an antisymmetric matrix for each $(x,t)$. A particular model for $\mathbb{S}$ in \cite{ShelRuz2013} can be written as
\begin{equation}
\left\{\begin{array}{rl}
\!\mathbb{S}&=\,\mu_1|\mathrm{B}_{s}|^{p-2}\mathrm{B}_{s}+\mu_2|\mathrm{R}|^{p-2}\mathrm{R}+\tau_\ast{\displaystyle\frac{\mathrm{B}_{\nu,p}}{|\mathrm{B}_{\nu,p}|}}, \quad \mbox{ if } \mathrm{B}_\nu\not=0\\ 
|\mathbb{S}|\! &\le\,\tau_\ast,  \quad \mbox{ if } \mathrm{B}_\nu=0,
\end{array}\right.
\label{SR explicit}
\end{equation}
where $\mu_1>0$, $\mu_2>0$ are (constant) viscosity coefficients, $
{\nu}\ge 0$ ({\em plug parameter}), 
\begin{equation}
\mathrm{B}_{\nu,p}=|\mathrm{B}_s|^{p-2}\mathrm{B}_s+\nu|\mathrm{R}|^{p-2}\mathrm{R} \qquad \mbox{and} \qquad \mathrm{B}_\nu=\mathrm{B}_{\nu,2}=\mathrm{B}_s+\nu\mathrm{R}.
\label{Bnu}
\end{equation}
In section \ref{mSR} we provide some additional details about the models \eqref{SR explicit}.

\smallskip

In this paper we are interested mostly in the Navier-Stokes system \eqref{NS} and  we consider $\Omega$ as a known function. That is, let us consider ${\Omega}$ as a given matrix function of $(x,t)$.

\smallskip

The first part of $\mathbb{S}$ in \eqref{SR explicit}, i.e., $\mu_1|\mathrm{B}_{s}|^{p-2}\mathrm{B}_{s}+\mu_2|\mathrm{R}|^{p-2}\mathrm{R}$, is the gradient of the function
\begin{equation}
U(X):=\frac{\mu _{1}}{p}|X_{s}|^{p}+\frac{\mu _{2}}{p}|X_{a}-\Omega|^{p}, \quad \mathrm{X}\in\mathbb{R}^{d\times d},
\label{U}
\end{equation}
at each matrix $\mathrm{B}$, for any $(x,t)$. However, the part $\frac{\mathrm{B}_{\nu,p}}{|\mathrm{B}_{\nu,p}|}$ has a potential (it is a gradient function, in the domain $\mathbb{R}^{d\times d}/\{\mathrm{X};\mathrm{X}_s+{\nu}(\mathrm{X}_a-\Omega)\not=\mathrm{0}\}$) if, and only if, $p=2$ and ${\nu}=0$ or $1$.\footnote{This fact is in accordance with Lemma 2.11 in the paper \cite{smr}, which states that the {\em plastic operator} $P(\mathrm{B}_s,\mathrm{R}):={\mathrm{B}_\nu}/{|\mathrm{B}_\nu|}$ is monotone if and only if ${\nu}=0$ or $1$. Further comments on this can be found in Section \ref{mSR}.} Indeed, this are necessary and sufficient conditions for the vector field 
$$
\frac{|\mathrm{X}_s|^{p-2}\mathrm{X}_s+{\nu}|\mathrm{X}_a+\mathrm{C}|^{p-2}(\mathrm{X}_a+\mathrm{C})}{\left|\, |\mathrm{X}_s|^{p-2}\mathrm{X}_s+{\nu}|\mathrm{X}_a+\mathrm{C}|^{p-2}(\mathrm{X}_a+\mathrm{C})\, \right|},
$$
where $\mathrm{C}\in\mathbb{R}^{d\times d}$ is fixed, to be a closed vector field in the space $\mathbb{R}^{d\times d}$ of $d\times d$ real matrices, seen $\mathbb{R}^{d\times d}$ as the direct sum of the subspaces of symmetric and antisymmetric matrices. (Compare with the vector field
$$
\frac{(|x|^{p-2}x,{\nu}|y+c|^{p-2}(y+c))}{\left|(|x|^{p-2}x,{\nu}|y+c|^{p-2}(y+c))\right|}
$$
in $\mathbb{R}^2$, where $c$ is a constant. It is easy to verify that this vector field is a closed vector field, in the domain $\mathbb{R}^2/\{(x,y); (x,{\nu} (y+c))\not=(0,0)\}$, iff $p=2$ and ${\nu}=0$ or $1$.) 

Setting $p=2$ in the term $\frac{\mathrm{B}_{{\nu},p}}{|\mathrm{B}_{{\nu},p}|}$ (taking $p=2$ only in this term in \eqref{SR explicit}), we obtain from \eqref{SR explicit} the tensors $\mathbb{S}$ for ${\nu}=0$ and ${\nu}=1$ given respectively by
\begin{equation}
\left\{\begin{array}{rl}
\!\mathbb{S}&=\,\mu_1|\mathrm{B}_{s}|^{p-2}\mathrm{B}_{s}+\mu_2|\mathrm{R}|^{p-2}\mathrm{R}+\tau_\ast{\displaystyle\frac{\mathrm{B}_{s}}{|\mathrm{B}_{s}|}}, \quad \mbox{ if } \mathrm{B}_s\not=0\\ 
|\mathbb{S}|\! &\le\,\tau_\ast,  \quad \mbox{ if } \mathrm{B}_s=0
\end{array}\right.
\label{SR0}
\end{equation}
and
\begin{equation}
\left\{\begin{array}{rl}
\!\mathbb{S}&=\,\mu_1|\mathrm{B}_{s}|^{p-2}\mathrm{B}_{s}+\mu_2|\mathrm{R}|^{p-2}\mathrm{R}+\tau_\ast{\displaystyle\frac{\mathrm{B}-\Omega}{|\mathrm{B}-\Omega|}}, \quad \mbox{ if } \mathrm{B}\not=\Omega\\ 
|\mathbb{S}|\! &\le\,\tau_\ast,  \quad \mbox{ if } \mathrm{B}=\Omega
\end{array}\right.
\label{SR1}
\end{equation}
(cf. (2.4)-(2.5), (2.20) and (2.30) in \cite{smr}).
These tensors also has a potential. More precisely, they are achieved respectively by the single condition $\mathbb{S}\in\partial V_i(\mathrm{B})$, $i=0,1$, where
$$
V_i(\mathrm{X})=U(\mathrm{X})+\tau_\ast W_i(\mathrm{X}),
$$
being $U$ the matrix function defined in \eqref{U},
$$
W_0(\mathrm{X})=|\mathrm{X}_s| \qquad \mbox{and} \qquad W_1(\mathrm{X})=|\mathrm{X}-\Omega|.
$$

We remark that the potentials $V_i$, $i=0,1$, as well as \eqref{U},  are functions depending on $(x,t)$ as well, since $\Omega$ does depend on $(x,t)$. To take into account stress tensors as these, we shall consider potentials $V\equiv V_{(x,t)}$ depending on $(x,t)$ as well. The condition $\mathbb{S}\in\partial V(\mathrm{B})$ for a potential $V\equiv V_{(x,t)}$ is understood in the pointwise sense with respect to almost all $(x,t)$, i.e.,  
\begin{equation}
\mathbb{S}(x,t)\in\partial V_{(x,t)}(\mathrm{B}(x,t)) \quad \mbox{for almost all} \quad (x,t).
\label{pointwise}
\end{equation}  

The symmetric and antisymmetric parts of $\mathbb{S}$ in \eqref{SR explicit} when $\mathrm{B}_\nu\not=0$ (equivalently, $\mathrm{B}_{\nu,p}\not=0$) are respectively equal to
$$
(\mu_1|\mathrm{B}_{s}|^{p-2}|\mathrm{B}_{{\nu},p}|+\tau_\ast)\frac{\, |\mathrm{B}_s|^{p-2}}{\!\! |\mathrm{B}_{{\nu},p}|}\mathrm{B}_{s}  
\quad \mbox{and} \quad
(\mu_2|\mathrm{R}|^{p-2}|\mathrm{B}_{{\nu},p}|+\tau_\ast)\frac{{\nu}|\mathrm{R}|^{p-2}}{|\mathrm{B}_{{\nu},p}|}\mathrm{R},
$$
then $|\mathbb{S}|$ equals the square root of
$$
\left[(\mu_1|\mathrm{B}_{s}|^{p-2}|\mathrm{B}_{{\nu},p}|+\tau_\ast)^2|\mathrm{B}_s|^{2(p-1)} + (\mu_2|\mathrm{R}|^{p-2}|\mathrm{B}_{{\nu},p}|+\tau_\ast)^2{\nu}^2|\mathrm{R}|^{2(p-1)}\right]\big/|\mathrm{B}_{{\nu},p}|^2
$$
and so is bigger than $\tau_\ast$. Thus, the plug region $\Omega^*_t$ for the model \eqref{SR explicit} is the set $\{x; \mathrm{B}_\nu=0\}$. In particular, for the models \eqref{SR0} and \eqref{SR1} ($\nu=0$ and $1$) we have $\Omega^*_t$ equals the sets $\{x; \mathrm{B}_s=0\}$ and $\{x; \mathrm{B}=\Omega\}$, respectively.

In this paper we consider the following generalization of the models \eqref{SR0} and \eqref{SR1}:
\begin{equation}
\left\{\begin{array}{l}
\mathbb{S}=
\mathrm{B}_{\mu ,p}+\widehat{\tau }_{\ast }\frac{\mathrm{B}_{{\nu} ,q}}{|\mathrm{B}^{\widehat{ \ }}_{{\nu} ,q}|^{\frac{2(q-1)}{q}}}, 
\quad \mbox{ if \ }\mathrm{B}_\nu\neq 0 \\ 
|\mathbb{S}|\leqslant \tau _{\ast }, 
\quad \mbox{ if \ }\mathrm{B}_\nu=0,
\end{array}\right.
\label{SR01g}
\end{equation}
with exponents $p,q\ge 2$, 
\begin{equation}
\widehat{\tau}_{\ast}=\frac{\tau_{\ast}}{\max (1,{\sqrt[q]{\nu}})}, \quad \tau_\ast,\, {\nu}\, \ge0,
\label{tau star hat}
\end{equation}
\begin{equation}
\mathrm{B}_{\mu ,p}=\mu _{1}|\mathrm{B}_{s}|^{p-2}\mathrm{B}_{s}+\mu _{2}|\mathrm{R}|^{p-2}\mathrm{R},
\label{B mu p}
\end{equation}
$\mu=(\mu_1,\mu_2), \quad \mu_1>0, \quad \mu_2\ge0$,
\begin{equation}
\mathrm{B}_{{\nu},q}=|\mathrm{B}_{s}|^{q-2}\mathrm{B}_{s}+{\nu} |\mathrm{R}|^{q-2}\mathrm{R}
\label{B nu q}
\end{equation}
and
\begin{equation}
\mathrm{B}^{\widehat{ \ }}_{{\nu},q}=|\mathrm{B}_{s}|^{\frac{q-2}{2}}\mathrm{B}_{s}+\sqrt{{\nu}} |\mathrm{R}|^{\frac{q-2}{2}}\mathrm{R}.
\label{B hat nu q}
\end{equation}
Our motivation for studying these tensors given by \eqref{SR01g} is purely a mathematical motivation. The point is that the modification in the tensors \eqref{SR explicit} through the expression  
$$
\widehat{\tau }_{\ast }\frac{\mathrm{B}_{{\nu},q}}{|\mathrm{B}^{\widehat{ \ }}_{{\nu} ,q}|^{\frac{2(q-1)}{q}}}
$$
results in tensors with a (convex) potential $V$, for any values of the above parameters, except that in the case $\nu=0$ we must assume also $\mu_2=0$ (in particular, for any ${\nu}>0$). More precisely, in Section \ref{OneDim} - see Corollary \eqref{S} - we show that the tensors given by \eqref{SR01g} can be achieved by the condition $\mathbb{S}(x,t)\in\partial V(\mathrm{B}(x,t))$, if $\nu>0$ or $\nu=\mu_2=0$, for the potential
\begin{equation}
V\equiv V_{(x,t)}=U(X)+\widehat{\tau }_{\ast }W(X)
\label{pot}
\end{equation}
with
\begin{equation}
\begin{array}{rl}
U(X)&\!\!=\,\frac{\mu _{1}}{p}|X_{s}|^{p}+\frac{\mu _{2}}{p}|\mathrm{R}|^{p}, \\
W(X)&\!\!=\,\left||\mathrm{X}_{s}|^{\frac{q-2}{2}}\mathrm{X}_{s}+\sqrt{{\nu}}|\mathrm{R}|^{\frac{q-2}{2}}\mathrm{R}\right|^{\frac{2}{q}}\\
    &\!\!=\,\sqrt[q]{|\mathrm{X}_{s}|^q+\nu |\mathrm{R}|^q},
\end{array}
\label{potential}
\end{equation}
where $\mathrm{R}=X_a-\Omega$.  

We do not have any physical or numerical experimentation to justify the Navier-Stokes systems \eqref{NS} with tensors \eqref{SR01g}. However, we believe that they are interesting from the mathematical point of view. In particular, in Section \ref{mSR}, we give a quite detailed description of the subdifferential of $V$ - see Theorem \ref{convexe1}.

More specifically, we solve the following three dimensional\footnote{Although our ibvp is set in three dimensional domain, many of our considerations in this paper can be formulated in the arbitrary dimension $d$. Thus, as much as the restriction on the dimension $d$ does not bring any big difficult we not restrict $d$ to a particular dimension.} initial boundary value problem (ibvp) (show the existence of a solution) - see Theorem \ref{theorem 1} below -, in a bounded domain $\mathcal{O}$ in $\mathbb{R}^{3}$\ with a boundary $\Gamma$ of class $C^{2}$, for any fixed time $T>0$, and with $\mathbb{S}$ given by \eqref{SR01g}: 
\begin{equation}
\left\{ 
\begin{array}{ll}
\mathbf{v}_{t}+\left( \mathbf{v\cdot \nabla }\right) \mathbf{v}=\mbox{div}\ 
\mathbb{T},\quad \mbox{div}\,\mathbf{v}=0 & \quad \text{in}\qquad \mathcal{O}
_{T}=\mathcal{O} \times (0,T) \\ 
\mathbf{v}\cdot \mathbf{n}=0,\;\quad (\mathbb{T}\,\mathbf{n}+\nabla g(\mathbf{v}))\cdot\bm{\tau}=0 & \quad \text{on}\qquad \Gamma _{T}=\Gamma \times
(0,T) \\ 
\mathbf{v}=\mathbf{v}_{0} & \quad \text{on}\qquad \mathcal{O} \times \{t=0\},%
\end{array}%
\right.  
\label{OneDeq1}
\end{equation}
where $\mathbf{n}$ is the vector field normal to $\Gamma$ pointing to the outside of $\Omega$ and $\bm{\tau}=\bm{\tau}({x})$ denotes any vector tangent to $\Gamma$ at the point ${x}\in\Gamma$. 

As for the boundary condition $(\mathbb{T}\,\mathbf{n}+\nabla g(\mathbf{v}))\cdot\bm{\tau}=0$, \ $g=g_{(x,t)}(\mathbf{v})$ is a Caratheodory function $g:\Gamma_T\times\mathbb{R}^3\to\mathbb{R}$, being differentiable and convex with respect to $\mathbf{v}$. It means that the tangential parts to $\Gamma$ of $\mathbb{T}\,\mathbf{n}$ and $\nabla g(\mathbf{v})=\nabla g_{(x,t)}(\mathbf{v}(x,t))$ coincide at any point $(x,t)\in\Gamma_T$.\footnote{The gradient operator $\nabla$ here acts only on the variable $\mathbf{v}$.} We also assume that $g$ satisfies the non-negativity and boundedness conditions:
\begin{equation}
0\le \nabla g_{(x,t)}(\mathbf{v})\cdot \mathbf{v}, \qquad |\nabla g_{(x,t)}(\mathbf{v})| \le \mathrm{c}|\mathbf{v}|,
\label{condition for g}
\end{equation}
for all $\mathbf{v}\in\mathbb{R}^3$ and $(x,t)\in\Gamma_T$, where $\mathrm{c}$ is a positive constant, independent of $\mathbf{v}$ and $(x,t)$. \ 
The particular case $g=\frac{1}{2}\alpha|\mathbf{v}|^2$ 
with $\alpha:\Gamma_T\to\mathbb{R}$ being a non-negative function in $L^\infty(\Gamma_T)$, satisfies \eqref{condition for g} and yields the Navier (friction/slip) boundary condition $(\mathbb{T}\,\mathbf{n}+\alpha \mathbf{v})\cdot\bm{\tau}=0$. 

The initial data $\mathbf{v}_{0}$ is a given function in the space $H$ defined below.

\

For the majority of this paper, we will, for simplicity, omit the explicit dependence of the potentials and the {\lq\lq}boundary function{\rq\rq} $g$ on the variables $x$ and $t$. That is, we will use the notation $V$ or $V(\mathrm{B})$ and $g(\mathbf{v})$ instead of $V_{(x,t)}(\mathrm{B})$ and $g_{(x,t)}(\mathbf{v})$.

We shall also use throughout the paper the following function spaces: 
\begin{equation}
\begin{array}{l}
H=\{\mathbf{v}\in L^2(\mathcal{O} ):\,\mbox{div }\mathbf{v}=0\;\ \text{ in}\;%
\mathcal{D}^{\prime }(\mathcal{O} ),\;\ \mathbf{v}\cdot \mathbf{n}=0\;\ \text{ in}%
\;H^{-1/2}(\Gamma )\}, \\ 
V=\{\mathbf{v}\in H^{1}(\mathcal{O} ):\,\mbox{div }\mathbf{v}=0\;\ \text{ a.e. in
\ }\mathcal{O} ,\;\ \mathbf{v}\cdot \mathbf{n}=0\;\ \text{ in}\;H^{1/2}(\Gamma
)\}, \\ 
V_{p}=\{\mathbf{v}\in V:\,\;\;|\nabla \mathbf{v}|\in L^p(\mathcal{O} )\}.%
\end{array}
\label{w}
\end{equation}%
where we have used standard notations of Lebesgue and Sobolev spaces and of
the space of regular distributions over $\mathcal{O} $. The space $V_{p}$ is
endowed with the norm 
\begin{equation*}
\Vert \mathbf{v}\Vert _{V_{p}}=\Vert \mathbf{v}\Vert _{L^{2}(\mathcal{O} )}+\Vert
\nabla \mathbf{v}\Vert _{L^{p}(\mathcal{O} )}.
\end{equation*}

\bigskip

We solve the ibvp \eqref{OneDeq1} in the sense of the following theorem. 
\begin{theorem}
\label{theorem 1} 
Let $V$ be the potential (function) given in \eqref{pot}-\eqref{potential}, being $\Omega\equiv\Omega(x,t)$ a matrix function in $L^p(\mathcal{O}_T;\mathbb{R}^{d\times d})$, and $g$ as above (a Caratheodory function from $\Gamma_T\times\mathbb{R}^3$ to $\mathbb{R}$, differentiable and convex with respect to the second variable and satisfying \eqref{condition for g}). Then for any $\mathbf{v}_{0}\in H$, $p\geqslant 2.2$ and $q\geqslant 2$, there exists a pair $(\mathbf{v},\mathbb{S)}$ of functions with $\mathbb{S}(x,t)\in\partial V_{(x,t)}(B(x,t))$ for almost all $(x,t)\in\mathcal{O}_T$ (in particular\footnote{See Corollary \ref{S}}, $\mathbb{S}$ satisfies \eqref{SR01g} for almost all $(x,t)\in\mathcal{O}_T$), such that 
\begin{equation}
\begin{array}{c}
\mathbf{v} \in L^{\infty }(0,T;H)\cap L^{p}(0,T;V_{p}),\quad \mathbf{v}
_{t}\in L^{p^{\prime }}(0,T;V_{p}^{\ast })  \qquad \mbox{and}
\\
\mathbb{S} \in L^{p^{\prime }}(\mathcal{O} _{T}),
\end{array} 
\label{regvw} 
\end{equation}
where $p^{\prime }=\frac{p}{p-1}$ and $V_{p}^{\ast }$ is the dual space of $%
V_{p}$. The pair $(\mathbf{v,}\mathbb{S)}$ satisfies the integral equality 
\begin{equation}
\begin{array}{c}
{\displaystyle\int }_{\!\!\!\mathcal{O} _{T}}[\mathbf{v}\cdot \partial _{t}{%
\boldsymbol{\varphi }}+\left( \mathbf{v\otimes v}-\mathbb{S}\right) :\nabla {%
\boldsymbol{\varphi }}]\,d{x}dt+{\displaystyle\int }_{\!\!\!\mathcal{O} }%
\mathbf{v}_{0}\cdot {\boldsymbol{\varphi }}({x},0)\,d{x} \\ 
={\displaystyle\int }_{\!\!\!\Gamma _{T}}\nabla g (\mathbf{v})\cdot {%
\boldsymbol{\varphi }}\,d\Gamma dt%
\end{array}
\label{weakForm}
\end{equation}%
for any function $\boldsymbol{\varphi }\in C^{1}(\overline{\mathcal{O} _{T}})$
such that $\boldsymbol{\varphi }(\cdot ,T)=0$,  $(\boldsymbol{\varphi }%
\cdot \mathbf{n})|\Gamma _{T}=0$ and $\mbox{div}\,\boldsymbol{\varphi }(\cdot ,t)=0$, for each $t\in [0,T]$.
\end{theorem}

\bigskip

Formally (for a smooth vector field $\mathbf{v}\in C^1({\overline{\mathcal{O}_T}})$), the Navier-Stokes system \eqref{NS} for $(x,t)\in{\mathcal{O}_T}$, $\mathbf{v}_{0}:=\mathbf{v}(\cdot,0)$, $\mbox{div}\,\mathbf{v}=0$ in $\mathcal{O}$ for any $t\in (0,T)$ and $\mathbf{v}\cdot \mathbf{n}=0$ on $\Gamma$ for any $t\in (0,T)$, is equivalent to the integral equality
\begin{equation}
\begin{array}{c}
{\displaystyle\int }_{\!\!\!\mathcal{O} _{T}}[\mathbf{v}\cdot \partial _{t}{
\boldsymbol{\varphi }}+\left( \mathbf{v\otimes v}-\mathbb{S}\right):\nabla{
\boldsymbol{\varphi }}]\,d{x}dt+{\displaystyle\int }_{\!\!\!\mathcal{O} }
\mathbf{v}_{0}\cdot {\boldsymbol{\varphi }}({x},0)\,d{x} \\ 
= -\ {\displaystyle\int }_{\!\!\!\Gamma _{T}}(\mathbb{T}\mathbf{n})\cdot {
\boldsymbol{\varphi }}\,d\Gamma dt,
\end{array}
\label{weakForm1}
\end{equation} 
where $\boldsymbol{\varphi }$ varies on the set of all vector functions in $C^{1}(\overline{\mathcal{O} _{T}};\mathbb{R}^n)$ vanishing at $t=T$ and divergence free in $\mathcal{O}$ for any $t\in (0,T)$ -- which yields \eqref{weakForm} by the boundary condition $(\mathbb{T}\,\mathbf{n}+\nabla g(\mathbf{v}))\cdot\bm{\tau}=0$ and the additional condition $(\boldsymbol{\varphi}\cdot\mathbf{n})|\Gamma_T=0$. Similarly, multiplying \eqref{NS}${}_1$ by $\mathbf{v}$ and integrating by parts we obtain (formally/rigorous for all divergence free $\mathbf{v}$ of class $C^1(\overline{\mathcal{O}_T})$) the {\lq}{\lq}energy equation{\rq\rq}
\begin{equation}
\frac{1}{2}\int_{\mathcal{O}}|\mathbf{v}(x,t)|^{2}d{x}
+\int_{\mathcal{O}_t}\mathbb{S}:{\nabla }\mathbf{v}d{x}ds-\int_{\Gamma_t}
(\mathbb{T}\mathbf{n})\cdot \mathbf{v}\,d{\Gamma } ds=\frac{1}{2}
\int_{\mathcal{O}}|\mathbf{v}_{0}|^{2}d{x},  
\label{energy eq}
\end{equation}
for any $t\in (0,T]$, where $\mathcal{O}_t:=\mathcal{O}\times (0,t)$ and $\Gamma_t:=\Gamma\times (0,t)$. \ The equations \eqref{weakForm1} and \eqref{energy eq} are essential elements for the analytic formulation and solvability of initial boundary values problems for the Navier-Stokes system. Particularly, we observe that when $\mathbb{S}=\nabla V(\mathrm{B})$, $\mathrm{B}=\nabla\mathbf{v}$, for some differentiable potential (function) $V$ (which is the case in the examples given above if we take $\tau_\ast=0$), the terms $\mathbb{S}:\nabla\boldsymbol{\varphi}$ and $\mathbb{S}:{\nabla}\mathbf{v}$ can be written in terms of the  directional derivatives $V'(\mathrm{X};\mathrm{Y})\equiv\partial V(X)/\partial Y=\nabla V(X):Y$ of $V$, i.e.,
$$
\mathbb{S}:\nabla\boldsymbol{\varphi}=V'(\mathrm{B};\nabla\boldsymbol{\varphi}) \qquad \mbox{and} \qquad \mathbb{S}:{\nabla}\mathbf{v}=V'(\mathrm{B};\mathrm{B}).
$$
When V is only convex, we could use one-sided directional derivatives here instead.\footnote{See section \ref{OneDim} for the definition and existence of one-sided directional derivatives for convex functions.} However, in solving ibvp's the velocity field is obtained by taking the limit of appropriate approximations and this arguably requires the linearity of the one-sided directional derivative $V'(X,Y)$ with respect to the {\lq}{\lq}direction{\rq}{\rq} $Y$. Thus, to possibly avoid extra technical difficulties, we may approximate our potential by differentiable potentials, but still preserving the good estimate satisfied by the potential. This is our general plan to solve an ibvp for the Navier-Stokes system having the viscous stress tensor in the subgradient of a convex function.

For concreteness,\footnote{A work in progress is underway to develop the plan outlined above into an arbitrary convex potential that satisfies appropriate conditions.} we have the following lemma regarding the potential~\eqref{pot}. 

\begin{lemma}
\label{apr copy(1)}
Let $\Omega$ be any fixed $d\times d$ antisymmetric real matrix. For given $p,q\geqslant 2$, let us consider the convex potential \eqref{pot}-\eqref{potential} and the approximations 
\begin{equation*}
\begin{array}{rll}
V^{n}(X) \ =&U(X) \ +&\!\!\!\!\!\!\!\!\!\!\!\!\!\!\!\!\!\!\!\!\!\!\!\!\!\!\!\!\!\!\!\!\!\!\!\!\!\!\!\!\!\!\!\!\!\!\!\!\!\widehat{\tau }_{\ast }W^n(X),\\
&U(X)=\frac{\mu _{1}}{p}|X_{s}|^{p}+\frac{\mu _{2}}{p}|\mathrm{R}|^{p} &\quad \mbox{\em (see \eqref{potential})},\\
&W^n(X) \ =&\!\!\!\!\!\!\!\!\!\!\!\!\!\!\!\!\!\!\!\!\!\!\!\!\!\!\!\!\!\!\!\!\!\!\!\!\!\!\!\!\!\!\sqrt[q]{|\mathrm{X}_{\nu,q}^{\widehat{ \ }}|^{2}+n^{-1}}\ =\sqrt[q]{|\mathrm{X}_{s}|^q+\nu |\mathrm{R}|^q+n^{-1}}
\end{array}
\end{equation*}
\begin{equation}
\mathbb{S}^{n}=\mathrm{X}_{\mu,p}+\widehat{\tau }_{\ast }\frac{\mathrm{X}
_{{\nu},q}}{\sqrt[q]{(|\mathrm{X}_{{\nu},q}^{\widehat{ \ }}|^{2}+n^{-1})^{q-1}}},
\label{approximate stress}
\end{equation}
$X\in \mathbb{R}^{d\times d}$, $n\in\mathbb{N}$. Here, \ $\mathrm{R}:=\mathrm{X}_{a}-{\Omega}$ (cf. \eqref{R}) and $\mathrm{X}_{\mu,p}$, $\mathrm{X}_{{\nu}^{2},q}$, $\mathrm{X}_{{\nu} ,q}^{\widehat{ \ }}$ are defined by \eqref{B mu p}, \eqref{B nu q}, \eqref{B hat nu q}, replacing $\mathrm{B}$ by $\mathrm{X}$. Then the following statements are true.

\begin{itemize}
\item[(a)] Both the potentials $V$ and $V^n$ satisfy the estimates
\begin{equation}
\mu_1 |X_s|^p - \mu_2 2^{p-2}|\Omega|^p - \tau_\ast |\Omega| \ \le \ V'(X;X)  \ \le \ c_1|X|^p + c_2|\Omega|^p + \tau_\ast |X|
\label{estimate for V and Vn}
\end{equation}
for any matrix $X\in\mathbb{R}^d\times\mathbb{R}^d$, for any $n\in\mathbb{N}$ (we can replace $V$ by $V^n$ in \eqref{estimate for V and Vn}), where 

\centerline{$c_1=\mu_1+2^{p-2}\mu_2(1+\frac{1}{p}) \quad \mbox{and} \quad c_2=2^{p-2}\mu_2 (1-\frac{1}{p})$.}

\item[(b)] For any\ given $X\in \mathbb{R}^{d\times d}$, we have that $\mathbb{S}^{n}=
\nabla V^{n}(\mathrm{X})$. Since $V^n$ is differentiable and convex, this statement is equivalent to the variational inequality $V^{n}(Y)-V^{n}(\mathrm{X})\geqslant \mathbb{S}^{n}:(Y-\mathrm{X})$, for all $Y\in \mathbb{R}^{d\times d}$. Moreover, 
\begin{equation}
\qquad\qquad |\mathbb{S}^n|\le \mu_1|X_s|^{p-1}+\mu_2|R|^{p-1}+\tau_\ast, \qquad \forall\, n, \ \forall\, X\in \mathbb{R}^{d\times d}.
\label{est Sn}
\end{equation}
\end{itemize}
\end{lemma}

The proof of this lemma is obtained by straightforward computations, some estimations, Theorem \ref{convexe1} and Corollary \ref{est DW}. We give its proof in Section \ref{Existence}. That section is devoted to solve the following approximation to the ibvp \eqref{OneDeq1}, where $\mathbb{T}^{n}=-p^{n}\,\mathrm{I}+\mathbb{S}^{n}$:
\begin{equation}
\left\{ 
\begin{array}{c}
\mathbf{v}_{t}^{n}+\left( \mathbf{v}^{n}\mathbf{\cdot \nabla }\right) 
\mathbf{v}^{n}=\mbox{div}\ \mathbb{T}^{n},\qquad \mbox{div}\,\mathbf{v}
^{n}=0,\qquad \text{in}\quad \mathcal{O} _{T}, \\ 
\mathbf{v}^{n}\cdot \mathbf{n}=0,\;\quad (\mathbb{T}^n\,\mathbf{n}+\nabla g(\mathbf{v}^n))\cdot\bm{\tau}=0\qquad \text{on }\ \Gamma
_{T}, \\ 
\mathbf{v}^{n}|_{t=0}=\mathbf{v}_{0}^{n}\qquad \text{in}\quad \mathcal{O}.
\end{array}
\right.
\label{approximate ibvp}
\end{equation}

Let us now sketch the proof of Theorem \ref{theorem 1}. The details are given in sections \ref{Existence} and \ref{Limit}. Once we obtain the sequence $(\mathbf{v}_n)$, the energy equation \eqref{energy eq}, the non-negativity condition in \eqref{condition for g} for the boundary function $g$, the coercivity condition (estimate from below) in \eqref{estimate for V and Vn} for the approximate potential $V^n$ and an appropriate version of Korn's inequatility - see Proposition \ref{Korn} - yield the estimates
\begin{equation}
\begin{array}{c}
\|\mathbf{v}^{n}\|_{L^{\infty}(0,T;H)}+\|\mathbf{v}^{n}\|_{L^{p}(0,T;V_{p})}\leqslant C,\\
\Vert\mathbb{S}^{n}\Vert_{L^{p^{\prime}}(\mathcal{O}_{T})}\leqslant C,\\
\Vert\partial _{t}\mathbf{v}^{n}\Vert_{L^{p^{\prime}}(0,T;\,W^{\ast})}\leqslant C,
\end{array}
\label{estimates}
\end{equation}
where $C$ denote positive constants that do not depend on $n$ (but may depend on $
p$, $q$, $\mathbf{v}_{0}$, $\mu _{1}$, ${\nu}$ and $\alpha$) and $W^{\ast}$ is dual space of\footnote{In Section \ref{Existence} we solve \eqref{approximate ibvp} using the Galerkin method. For simplicity of calculations, we choose a basis of appropriate eigenfunctions in the space $W$.} $W:=W^{3,2}(\mathcal{O} )\cap V$.
Then, we have the following convergences (up to subsequences), with the help of the continuity of the trace map from $V_p\equiv W^{1,p}(\mathcal{O})$ to $L^p(\Gamma)$, the Aubin-Lions-Simon compactness result (Lemma \ref{ALS}) and the boundedness condition in \eqref{condition for g} for the gradient of the boundary function $g$:
\begin{equation}
\begin{array}{clll}
\mathbf{v}^{n} &\rightharpoonup &\mathbf{v}\quad \mbox{ weakly-$\ast$ in }
& L^{\infty }(0,T;H),\\
\mathbf{v}^{n} &\rightharpoonup &\mathbf{v}\quad \mbox{ weakly in }
& L^{p}(0,T;V_{p})\cap L^p(\Gamma _{T}),\\
\mathbb{S}^{n} &\rightharpoonup &\mathbb{S}\quad \mbox{ weakly in}
& L^{p^{\prime }}(\mathcal{O} _{T}),\\
\nabla g(\mathbf{v}^{n}) &\rightharpoonup &\mathbf{g}\quad \mbox{ weakly in}
& L^p(\Gamma_{T}),  
\end{array}
\label{conv introducao}
\end{equation}
for some $\mathbb{S}\in L^{p^{\prime }}(\mathcal{O} _{T})\equiv L^{p^{\prime }}(\mathcal{O} _{T}; \mathbb{R}^{3\times 3})$, $\mathbf{g}\in L^p(\Gamma_{T})\equiv L^p(\Gamma_{T};\mathbb{R}^3)$, and
\begin{equation}
\mathbf{v}^{n}\rightarrow \mathbf{v}\quad \mbox{strongly in }~L^{2}(\mathcal{O} _{T}).
\label{convergence introducao}
\end{equation}
Moreover, 
\begin{equation}
\partial _{t}\mathbf{v}\in L^{p^{\prime }}(0,T;W^{\ast }).  
\label{aa introducao}
\end{equation}
Thus, we can take $n\to\infty$ in the weak formulation of \eqref{approximate ibvp} - see \eqref{vnn} - and obtain \eqref{weakForm} with $\mathbf{g}$ in place of $\nabla g(\mathbf{v})$. At this point, it will remain to show only that $\mathbb{S}\in\partial V(\mathrm{B})$ and $\mathbf{g}=\nabla g(\mathbf{v})$. We can show this by a standard argument, which we explain in the sequel for the convenience of the reader.
Although the argument is standard, we would like to draw attention to the fact that here enters the condition $p\ge 2.2=11/5$ (i.e., in this part of the proof of Theorem \ref{theorem 1} we utilize Lemma \ref{1 copy(1)}).

In fact, the main ingredients to show that $\mathbb{S}\in\partial V(\mathrm{B})$ and $\mathbf{g}=\nabla g(\mathbf{v})$ appears to be (not necessarily in the following order) the lower semicontinuity of functionals of the form {\lq\lq}$\mathcal{F}(u)=\int F(x,u(x))dx${\rq\rq} in the weak topology \cite[Theorem 4.3, p. 123]{giusti}, the condition $p\ge 2.2$ and the energy equation \eqref{energy eq} for $\mathbf{v}^n$, which can be written as
\begin{equation}
\frac{1}{2}\int\limits_{\mathcal{O}} |\mathbf{v}^{n}(t)|^{2}-|\mathbf{v}
_{0}^{n}|^{2} \,d{x}=-\int\limits_{\mathcal{O} _{t}} \mathbb{S}^{n}:\mathrm{B}^{n} \,\,d{x}\,ds - \int\limits_{\Gamma_t}\nabla g(\mathbf{v}^{n})\cdot\mathbf{v}^n\,d{\Gamma }\,ds,  
\label{s2 introducao}
\end{equation}
for all $t\in (O,T)$. Indeed, from \cite[Theorem 4.3, p. 123]{giusti} and the convexity of the functions $V$, $g$, we have for any $X\in L^2(\mathcal{O}_T)$, ${y}\in L^2(\Gamma_T)$ and $t\in (0,T)$:
\begin{eqnarray}
&\int\limits_{\mathcal{O}_{t}}V(X)-V(\mathrm{B})\,d{x}ds 
+ \int\limits_{\Gamma_{t}}g({y})-g(\mathbf{v})\,d\Gamma ds  \notag\\ 
&\geqslant \liminf_{n\rightarrow \infty} [\int\limits_{\mathcal{O}_{t}}V^{n}(X)-V^{n}(\mathrm{B}^{n})\,d{x}ds + \int\limits_{\Gamma_{t}} g({y})-g(\mathbf{v}^n)\,d\Gamma ds] \notag \\
&\geqslant \liminf_{n\rightarrow\infty} [\int\limits_{\mathcal{O}_{t}}\mathbb{S}^{n}:(X-\mathrm{B}^n)\,d{x}ds + \int\limits_{\Gamma_{t}} \nabla g(\mathbf{v}^n)\cdot ({y}-\mathbf{v}^n)\,d\Gamma ds] \notag \\
&= \int\limits_{\mathcal{O} _{t}}\mathbb{S}:X \, dxds + \int\limits_{\Gamma_{t}}  \mathbf{g}\cdot {y} d\Gamma dt \notag \qquad\qquad\qquad\qquad\qquad\qquad\qquad\quad\\ 
& \qquad\qquad + \liminf_{n\rightarrow\infty} [-\int\limits_{\mathcal{O} _{t}}\mathbb{S}^{n}:\mathrm{B}^n\,d{x}ds - \int\limits_{\Gamma_{t}} \nabla g(\mathbf{v}^n)\cdot \mathbf{v}^n\,d\Gamma ds];
\label{eq0.0 introducao}
\end{eqnarray}
substituting \eqref{s2 introducao} here, we obtain
\begin{eqnarray}
&\int\limits_{\mathcal{O}_{t}}V(X)-V(\mathrm{B})\,d{x}ds 
+ \int\limits_{\Gamma_{t}}g({y})-g(\mathbf{v})\,d\Gamma ds  \notag\\
\geqslant & \int\limits_{\mathcal{O} _{t}}\mathbb{S}:X \, dxds + \int\limits_{\Gamma_{t}}  \mathbf{g}\cdot {y} d\Gamma dt +\liminf_{n\rightarrow\infty}[\frac{1}{2}\int\limits_{\mathcal{O} } |\mathbf{v}^{n}(t)|^{2}-|\mathbf{v}_{0}^{n}|^{2} \,d{x}]\notag;
\end{eqnarray}
then, using the strong convergence \eqref{convergence introducao} and $\mathbf{v}^n_0\to\mathbf{v}_0$ also strongly, it follows that
\begin{equation}
\begin{array}{ll}
&\int\limits_{\mathcal{O}_{t}}V(X)-V(\mathrm{B})\,d{x}ds 
+ \int\limits_{\Gamma_{t}}g({y})-g(\mathbf{v})\,d\Gamma ds \\
\geqslant & \int\limits_{\mathcal{O} _{t}}\mathbb{S}:X \, dxds + \int\limits_{\Gamma_{t}}  \mathbf{g}\cdot {y} d\Gamma dt + \frac{1}{2}\int\limits_{\mathcal{O} } |\mathbf{v}(t)|^{2}-|\mathbf{v}_{0}|^{2} \,d{x},
\end{array}
\label{almost there}
\end{equation}
for almost all $t\in (0,T)$. We now apply the 
counterpart formula/energy equation \eqref{s2 introducao} for $\mathbf{v}$ in the form
\begin{equation}
{\textstyle
\frac{1}{2}\int\limits_{\mathcal{O} } |\mathbf{v}(t)|^{2}-|\mathbf{v}
_{0}|^{2} \,d{x}=-\int\limits_{\mathcal{O} _{t}} \mathbb{S}:\mathrm{B}\,\,d{x}\,ds - \int\limits_{\Gamma_t}\mathbf{g}\cdot\mathbf{v}\,d{\Gamma }\,ds.} 
\label{baby s2 introducao}
\end{equation}   
Eventually, this formula evolves into the same formula \eqref{s2 introducao}, with the $n$'s omitted, as $\mathbf{g}$ becomes $\nabla g(\mathbf{v})$. It can be obtained 
from the weak formulation \eqref{weakForm} with $\mathbf{g}$ in place of $\nabla g(\mathbf{v})$, by choosing test functions approaching $\mathbf{v}$. In fact, using Lemma \ref{1 copy(1)} we can put $\boldsymbol{\varphi}=\zeta_\varepsilon(t)\mathbf{v}$ with $\zeta_\varepsilon\stackrel{\varepsilon\to0+}{\rightharpoonup}\chi_{(0,t)}$ in the equation 
\begin{eqnarray}
&\int\limits_{\mathcal{O} _{T}}[ \mathbf{v}\partial _{t}{\boldsymbol{\varphi }}+(\mathbf{v\otimes v}-\mathbb{S}) :\nabla {\boldsymbol{\varphi }}]dxdt +\int\limits_{\mathcal{O} }\mathbf{v}_{0}{\boldsymbol{\varphi }}(0)dx =\int\limits_{\Gamma_T}\mathbf{g}\cdot {\boldsymbol{\varphi}}d\Gamma dt
\notag
\end{eqnarray}
to get
\begin{eqnarray}
&-\int\limits_{\mathcal{O}}|\mathbf{v}(t)|^2dx + \int_0^t\langle \mathbf{v}_t,\mathbf{v}\rangle_{V_p^*,V_p}ds
-\int\limits_{\mathcal{O} _{t}}\mathbb{S} :\nabla \mathbf{v}dxds
 +\int\limits_{\mathcal{O}}|\mathbf{v}_{0}|^2dx 
=\int\limits_{\Gamma_t}\mathbf{g}\cdot\mathbf{v}d\Gamma ds 
\notag
\end{eqnarray}
for almost all $t\in (0,T)$. But 
$\int\limits_{\mathcal{O}} |\mathbf{v}(t)|^{2}-|\mathbf{v}_0|^{2} \,d{x}=\int_0^t 2\langle \mathbf{v}_t,\mathbf{v}\rangle_{V_p^*,V_p}ds$, so
\begin{eqnarray}
&-\frac{1}{2}\int_{\mathcal{O}} |\mathbf{v}(t)|^{2}-|\mathbf{v}_0|^{2} \,d{x} -\int_{\mathcal{O} _{t}}\mathbb{S} :\nabla \mathbf{v}
 =\int\limits_{\Gamma_t}\mathbf{g}\cdot\mathbf{v}d\Gamma ds,
\notag
\end{eqnarray}
i.e., we have \eqref{baby s2 introducao}. From \eqref{almost there} and \eqref{baby s2 introducao} we have that
\begin{equation}
\begin{array}{ll}
&\int\limits_{\mathcal{O}_{t}}V(X)-V(\mathrm{B})\,d{x}ds 
+ \int\limits_{\Gamma_{t}}g({y})-g(\mathbf{v})\,d\Gamma ds \\
\geqslant & \int\limits_{\mathcal{O} _{t}} \mathbb{S}:(\mathrm{X}-\mathrm{B})\,\,d{x}\,ds + \int\limits_{\Gamma_t}\mathbf{g}\cdot({y}-\mathbf{v})\,d{\Gamma }\,ds.
\end{array}
\end{equation}
Since the matrix function $X\in L^{2}(\mathcal{O} _{T})$ and the vector function ${y}\in L^2(\Gamma_T)$ are arbitrary, we can choose in this inequality 
$X=\mathrm{B}+\varepsilon \mathrm{Z}$, ${y}=\mathbf{v}+\varepsilon \mathbf{z}$ \ for any positive $\varepsilon$, any matrix function $\mathrm{Z}\in L^{2}(\mathcal{O} _{T})$ any vector function $\mathbf{z}\in L^{2}(\Gamma_{T})$,  \ which gives 
\begin{equation*}
\begin{array}{rl}
           &\int\limits_{\mathcal{O}_{t}}\frac{V(\mathrm{B}+\varepsilon \mathrm{Z})-V(\mathrm{B})}{\varepsilon}\,d{x}ds + \int\limits_{\Gamma_{t}}\frac{g(\mathbf{v}+\varepsilon \mathbf{z})-g(\mathbf{v})}{\varepsilon}\,d{\Gamma}ds\\
\geqslant & \qquad\int\limits_{\mathcal{O}_{t}}\mathbb{S}:\mathrm{Z}\,d{x}dt + \int\limits_{\Gamma_{t}}\mathbf{g}:\mathbf{z}\,d{\Gamma}ds.
\end{array}
\end{equation*}
So, taking the limit when $\varepsilon\to 0+$ and passing the limit over the sign of integration, by the Vitali convergence theorem, we obtain the following inequality for the one-sided directional derivatives of the convex functions $V$ and $g$:
\begin{equation*}
{\textstyle\int\limits_{\mathcal{O}_{t}}V^{\prime}(\mathrm{B};\mathrm{Z})\,d{x}ds + \int\limits_{\Gamma_{t}}g^{\prime}(\mathbf{v};\mathbf{z})\,d{\Gamma}ds \geqslant \int\limits_{\mathcal{O}_{t}}\mathbb{S}:\mathrm{Z}\,d{x}ds + \int\limits_{\Gamma_{t}}\mathbf{g}:\mathbf{z}\,d{\Gamma}ds,}
\end{equation*}%
for any matrix function $\mathrm{Z}\in L^{2}(\mathcal{O}_{T})$ and any vector function $\mathbf{z}$. Now, using the positive homogeneity of one-sided derivatives of convex functions with respect to the {\lq\lq}directions{\rq\rq}, this inequality implies indeed in 
\begin{equation*}
{\textstyle\int\limits_{\mathcal{O}_{t}}V^{\prime}(\mathrm{B};\mathrm{Z})\xi \,d{x}ds + \int\limits_{\Gamma_{t}}g^{\prime}(\mathbf{g};\mathbf{z})\xi \,d{\Gamma}ds
\geqslant \int\limits_{\mathcal{O}_{t}}(\mathbb{S}:\mathrm{Z})\xi \,d{x}ds
+\int\limits_{\Gamma_{t}}(\mathbf{g}:\mathbf{z})\xi \,d{\Gamma}ds,}
\end{equation*}%
for any positive bounded mensurable function $\xi$. Therefore, $V^{\prime }(\mathrm{B};\mathrm{Z})\geqslant \mathbb{S}:\mathrm{Z}$ and $g^{\prime }(\mathbf{v};\mathbf{z})\geqslant \mathbf{g}\cdot \mathbf{z}$ \ for any matrix function $\mathrm{Z}\in L^{2}(\mathcal{O} _{T})$ and vector function $\mathbf{z}\in L^2(\Gamma_T)$, and thus, $\mathbb{S}(x,t)\in \partial V_{(x,t)}(B(x,t))$ and $\mathbf{g}(x,t)=\nabla g_{(x,t)}(\mathbf{v}(x,t))$, by Theorem \ref{dir}.

\

While our introduction has been extensive, we trust it serves as a valuable motivation and facilitates the reading of the paper. Notably, the last part above constitutes a detailed fraction of the proof of Theorem \ref{theorem 1}, thereby abbreviating Section \ref{Limit}.

Yet before closing it, we must comment on some previous results concerning 
the solvability of the asymmetric Herschel-Bulkley/Cosserat-Bingham fluid equations. 
In \cite{smr} it was obtained the solvability result of the {\bf{stationary}} boundary value problem for the Cosserat-Bingham model with the viscous stress tensor given by \eqref{SR1}, the homogeneous (null) Dirichlet boundary  condition and $p>6/5$. In \cite{ShelChem}, a {\bf non-stationary} boundary value problem for the asymmetric Herschel-Bulkley fluid with $p=q=2$ was studied in the one-dimensional domain, and, more recently, in \cite{Anderson-Nikolai} it was solved  the three dimensional case for the Cosserat-Bingham model \eqref{SR01g} with also $p=q=2$ and the homogeneous (null) Dirichlet boundary condition. Comparing \cite{Anderson-Nikolai} with our results in this paper, we loose in Theorem \ref{theorem 1} regarding the restriction $p\ge 2.2$, but this restriction simplifies the proof quite a bit, since in this case we can use Lemma \ref{1 copy(1)} and so we do not need to go over two levels of approximations as done there. More precisely, in \cite{Anderson-Nikolai} to prove their main result the authors take first the limit over Galerkin approximations for then take the limit over the approximations for the potential and viscous stress tensor. Here, we can take the same parameter for both approximations and take the limit at once. We believe it is worth to deliver our result with this restriction and show our proof for this case. In \cite{Anderson-Nikolai}, as well as in \cite{smr}, the full Cosserat-Bingham model is considered. Here, we fix a micro-rotational velocity tensor ${\Omega}$ and proceed the analysis allowing potentials/viscous stress tensors depending on the $x,t$ variables as well. This procedure may suggest an increment in the approach to solve the full Cosserat-Bingham fluid equations. Differently from \cite{smr} and \cite{Anderson-Nikolai}, our boundary condition is of Navier boundary condition type having a somewhat generic form. Overall, we exploit carefully the structure of the viscous stress tensor as the subgradient of a convex function and, as a byproduct, we give a quite complete characterization of the subgradient of the potential \eqref{pot}-\eqref{potential} - see Theorem \ref{convexe1}. Finally, we observe that our result covers any values of the exponents $p\ge 2.2$ and $q\ge2$.

To truly close this introduction, we provide a brief description of the organization of this paper. In section \ref{mSR} we give more details on the model \eqref{SR explicit} and the relation with \eqref{SR01g}. In section \ref{OneDim} we collect some technical results that are used to prove theorems \ref{theorem 1} and \ref{convexe1}. 
In section \ref{Existence} we solve the approximated problem \eqref{y1}. Also we derive a priori estimates for the solution of \eqref{y1}, which are independent on $n$. Finally, in section \ref{Limit} we prove Theorem \ref{theorem 1},
applying the Lions-Aubin compactness theorem, a priori estimates of 
Section \ref{Existence} and the theory of monotone operators.  

\section{The systems \eqref{SR explicit} and the modified version \eqref{SR01g}}
\label{mSR}

Let us consider {the {\lq\lq}asymmetric{\rq\rq} Herschel-Bulkley
fluid model or, more generally, the \emph{Cosserat-Bingham} fluid model,
proposed in \cite{ShelRuz2013}, where the rotation of particles is taken
into account, being described by a \emph{micro-rotational velocity tensor}
denoted by ${\Omega}$. In this model the tensor $\mathbb{S}$ is given
implicitly by the relation (3.8) in \cite{ShelRuz2013}, i.e. 
\begin{equation}  \label{CB}
\mu_1(a_1+|\mathrm{B}_{s}|)^{p-2}\left((|\mathbb{S}|-\tau_\ast)_{+}+\tau_%
\ast\right)\mathrm{B}_0=(|\mathbb{S}|-\tau_\ast)_{+}\mathbb{S},
\end{equation}
where $\mu_1, \tau_\ast>0$, $a_1\ge0$ and\footnote{Four our convenience in this paper, we have changed the notation {\lq\lq}$q${\rq\rq} in the exponents used in \cite{ShelRuz2013} to the letter {\lq\lq}$p${\rq\rq} and changed $2\mu_1$ to $\mu_1$.
} 
\begin{equation}  \label{CB cont'ed}
\mathrm{B}_0=\mathrm{B}_{s}+\varepsilon\mathrm{R}, \quad \mathrm{R}=\mathrm{B%
}_{a}-{\Omega}, \quad \varepsilon=\frac{\mu_2}{\mu_1}\frac{(a_2+|%
\mathrm{R}|)^{p-2}}{(a_1+|\mathrm{B}_{s}|)^{p-2}}, \quad \mu_2>0, \ a_2\ge0.
\end{equation}
Solving \eqref{CB}-\eqref{CB cont'ed}, first we observe that if $|\mathbb{S}%
|\le\tau_\ast$ then, from \eqref{CB}, it follows that 
\begin{equation*}
\mu_1(a_1+|\mathrm{B}_{s}|)^{p-2}\tau_\ast\mathrm{B}_0=0,
\end{equation*}
therefore, $\mathrm{B}_0=0$ and so, also $\mathrm{B}_{s}=\mathrm{R}=0$. (Here, we have
used the fact that ${\Omega}$ is skew-symmetric - see \cite{ShelRuz2013}.) \ On the other hand, if $|\mathbb{S}|>\tau_\ast$, also from 
\eqref{CB}, we have 
\begin{equation}  \label{CB cont''ed}
\mu_1(a_1+|\mathrm{B}_{s}|)^{p-2}|\mathbb{S}|\mathrm{B}_0=(|\mathbb{S}%
|-\tau_\ast)\mathbb{S},
\end{equation}
then $\mathrm{B}_0\not=0$ and, taking the norm on both sides of this
equation, we obtain that 
\begin{equation*}
\frac{|\mathbb{S}|-\tau_\ast}{\mu_1(a_1+|\mathrm{B}_{s}|)^{p-2}}=|\mathrm{B}%
_0|.
\end{equation*}
Substituting this relation in \eqref{CB cont''ed}, it follows that 
\begin{equation*}
\mathbb{S}=|\mathbb{S}|\frac{\mathrm{B}_0}{|\mathrm{B}_0|}=\tau_\ast\frac{%
\mathrm{B}_0}{|\mathrm{B}_0|} + \mu_1(a_1+|\mathrm{B}_{s}|)^{p-2}\mathrm{B}%
_0,
\end{equation*}
then, using \eqref{CB cont'ed}, we arrive at 
\begin{equation*}
\mathbb{S}=\mu_1(a_1+|\mathrm{B}_{s}|)^{p-2}\mathrm{B}_{s}+\mu_2(a_2+|%
\mathrm{R}|)^{p-2}\mathrm{R}+\tau_\ast\frac{\mathrm{B}_0}{|\mathrm{B}_0|}.
\end{equation*}
Summing up, we have that \eqref{CB}-\eqref{CB cont'ed} is equivalent to 
\begin{equation}  \label{CB explicit}
\begin{array}{rl}
|\mathbb{S}|>\tau_\ast \Leftrightarrow & \mathbb{S}=\mu_1(a_1+|\mathrm{B}%
_{s}|)^{p-2}\mathrm{B}_{s}+\mu_2(a_2+|\mathrm{R}|)^{p-2}\mathrm{R}+\tau_\ast{%
\displaystyle\frac{\mathrm{B}_0}{|\mathrm{B}_0|}}, \\ 
|\mathbb{S}|\le\tau_\ast \Leftrightarrow & \mathrm{B}_0=0.%
\end{array}
\end{equation}
The relations \eqref{CB explicit} are in \cite[(3.13)]{ShelRuz2013} for the case $p=2$. See also \cite{smr} - in particular, the equation (2.20), where is considered more general assumptions for the stress tensor $\mathbb{S}$. In particular, \eqref{CB explicit}  for $a_1=a_2=0$  yields 
\begin{equation}  \label{HB_a}
\begin{array}{rl}
|\mathbb{S}|>\tau_\ast \Leftrightarrow & \mathbb{S}=\mu_1|\mathrm{B}
_{s}|^{p-2}\mathrm{B}_{s}+\mu_2|\mathrm{R}|^{p-2}\mathrm{R} \\ 
& \ \ \ \ \ \ \ \ +\tau_\ast{\displaystyle\frac{|\mathrm{B}_{s}|^{p-2}\mathrm{B}_{s}+\frac{\mu_2}{\mu_1}|\mathrm{R}|^{p-2}\mathrm{R}}{\left||\mathrm{B}_{s}|^{p-2}\mathrm{B}_{s}+\frac{\mu_2}{\mu_1}|\mathrm{B}_{a}|^{p-2}\mathrm{R}\right|}} \\ 
|\mathbb{S}|\le\tau_\ast \Leftrightarrow & \mathrm{B}_{s}+\frac{\mu_2}{\mu_1}\mathrm{R}=0,
\end{array}
\end{equation}
which is the model \eqref{SR explicit}, if ${\nu}=\mu_2/\mu_1$. In this paper, as well as in \cite{smr}, we consider ${\nu}$ as a independent parameter of $\mu_2$ or $\mu_1$.

Comparing the tensor $\mathbb{S}$ in the model \eqref{SR explicit} with the tensor $\mathbb{S}$ in our system \eqref{SR01g}, we see that the difference between them is in the parts 
\begin{equation*}
{\tau }_{\ast }\frac{\mathrm{B}_{{\nu},p}}{|\mathrm{B}_{{\nu},p}|}={\tau}_{\ast}\frac{|\mathrm{B}_{s}|^{p-2}\mathrm{B}_{s}+{\nu}|\mathrm{B}_{a}|^{p-2}\mathrm{B}_{a}}{\sqrt{|\mathrm{B}_{s}|^{2(p-1)}+{\nu}^2|\mathrm{B}_{a}|^{2(p-1)}}},
\end{equation*}
in \eqref{SR explicit}, and 
\begin{equation*}
\widehat{\tau}_{\ast}\frac{\mathrm{B}_{{\nu},q}}{|\mathrm{B}_{{\nu},q}^{\widehat{ \ }}|^
\frac{2(q-1)}{q}}=\frac{\tau_\ast}{\max{(1,{\sqrt[q]{\nu}})}}\frac{|\mathrm{B}_{s}|^{q-2}\mathrm{B}_{s}+{\nu}|\mathrm{B}_{a}|^{q-2}\mathrm{B}_{a}}{(\sqrt[q]{|\mathrm{B}_{s}|^q+{\nu}|\mathrm{B}_{a}|^q})^{q-1}},
\end{equation*}
in \eqref{SR01g}. Observe that these two parts coincide when $q=p=2$
and ${\nu}=0,1$. For this reason, we can say that our family of systems 
\eqref{SR01g} contains a class of symmetric and asymmetric Herschel-Bulkley type equations. Notice that in the first part, $\mathrm{B}_{\mu,p}$, which is the same for 
\eqref{SR explicit} and \eqref{SR01g}, there is no restriction on the parameter $p$ in order to have the coincidence in the above parts when $q=p=2$ and ${\nu}=0,1$.}

We would like to stress that our motivation for considering the general tensor $\mathbb{S}$ satisfying \eqref{SR01g} is pure mathematical, i.e. although, when ${\nu}\not=0,1$ or $q>2$ the systems \eqref{SR01g} might not have any relation to some possibly real fluid, from a pure mathematical point of view we find it interesting to study the problem \eqref{OneDeq1} with this type of more general tensors. In particular, let us make the following point. Borrowing the notation and terminology used in \cite{smr}, consider the operators 
\begin{equation}
\mathbf{\mathrm{P}}_{{\nu},q}(\mathrm{B}_{s},\mathrm{R})={\displaystyle\frac{\mathrm{B}_{{\nu},q}}{|\mathrm{B}_{{\nu},q}|}}=\frac{|\mathrm{B}_{s}|^{q-2}\mathrm{B}_{s}+{\nu}|\mathrm{R}|^{q-2}\mathrm{R}}{\sqrt{|\mathrm{B}_{s}|^{2(q-1)}+{\nu}^2|\mathrm{R}}|^{2(q-1)}},
\label{po}
\end{equation}
\begin{equation}  \label{mpo}
\mathbf{\mathcal{P}}_{{\nu},q}(\mathrm{B}_{s},\mathrm{R})={\displaystyle\frac{\mathrm{B}_{{\nu},q}}{|\mathrm{B}_{{\nu},q}^{\widehat{ \ }}|^\frac{2(q-1)}{q}}}=\frac{|\mathrm{B}_{s}|^{q-2}\mathrm{B}_{s}+{\nu}|\mathrm{R}|^{q-2}\mathrm{R}}{(\sqrt[q]{|\mathrm{B}_{s}|^q+{\nu}|\mathrm{R}|^q})^{q-1}
}
\end{equation}
and call them \emph{plastic operator} and \emph{modified plastic operator},
respectively. With these notations, the Shelukhin-R\r{u}\v{z}i\v{c}ka
tensors \eqref{SR explicit} and ours \eqref{SR01g} can be written, respectively, as 
\begin{equation*}
\mathbb{S}=\mathrm{B}_{\mu,p}+\tau_\ast\mathbf{\mathrm{P}}_{{\nu},p}(\mathrm{B}_{s},\mathrm{R})
\end{equation*}
and 
\begin{equation*}
\mathbb{S}=\mathrm{B}_{\mu,p}+\widehat{\tau_\ast}\mathbf{\mathcal{P}}_{{\nu},q}(\mathrm{B}_{s},\mathrm{R}),
\end{equation*}
when $|\mathbb{S}|>\tau_\ast$. As we observed above, we have that ${\mathrm{P}}_{{\nu},q}=\mathbf{\mathcal{P}}_{{\nu},q}$ if ${\nu}=0,1$ and $q=2$. 
It was proved in \cite[Lemma 2.11]{smr} that the operator ${\mathrm{P}}_{{\nu},q}$, when $q=2$, is monotone if and only if ${\nu}=0,1$. In this paper, we shall prove that the \emph{modified plastic operator} $\mathbf{\mathcal{P}}_{{\nu},q}$ is monotone for any values of ${\nu}\ge0$ and $q\ge2$. See Corollary \ref{mono}.

Furthermore, each tensor in our family of tensors \eqref{SR01g} is a subgradient of the potential \eqref{pot}-\eqref{potential} (see Corollary \ref{S}), the subdifferential of which we characterize quite completely in Theorem \ref{convexe1}. 

\smallskip

In the proof of Theorem \ref{convexe1} we employ some basic facts from Convex Analysis, which we gather in the next section.

\section{Auxiliary results}
\label{OneDim}

In this section, we present some very useful facts for our paper.
We begin by recalling two basic theorems on convex analysis and making two remarks.

\begin{theorem}
\label{dir der} 
{\em (See \cite[Theorem 23.1]{R} or \cite[Proposition 17.2]{BC}.)}\\ 
Let $f:\mathbb{R}^m\to\mathbb{R}$ be a convex function. Then, for each ${y}\in \mathbb{R}^{m}$ there exists the one-sided directional derivative $f^{\prime }({x};{y})$ of $f$ at ${x}$ with respect to the vector~${y}$:
\begin{equation} 
f^{\prime }({x};{y}):=\lim_{\lambda\to 0+}\frac{f(x+\lambda y)-f(x)}{\lambda}=\inf_{\lambda>0}\frac{f(x+\lambda y)-f(x)}{\lambda}.
\label{def f'}
\end{equation}
\end{theorem}

The proof of this theorem consists in the elementary (however very important) observation that the difference quotient $(f({x}+\lambda {y})-f({x}))/\lambda $ is a non-decreasing function of $\lambda >0$, so that the above limit exists and coincides with the above infimum.

\begin{remark}
\label{hom} For any function $f:\mathbb{R}^{m}\rightarrow \mathbb{R}$ positively homogeneous of order 1, we have $f^{\prime }(0;{y})=f({y})$, forall $y\in \mathbb{R}^{m}$. {\em Indeed, in this case one has $f^{\prime }(0;y)=\lim_{\lambda\to 0+}\lambda^{-1}\left(f(\lambda y)-f(0)\right)=\lim_{\lambda\to 0+}\lambda^{-1}\left(\lambda f(y)\right)=f(y)$.}
\end{remark}

\medskip

Let us remember the concept of subdifferential.

\begin{definition}
\label{subd} 
A vector ${x}^{\ast }$ is said to be a subgradient of $f $ at ${x}$ if 
\begin{equation*}
f({y})\geqslant f({x})+{x}^{\ast }\cdot({y}-%
{x}),\quad \forall {y}\in \mathbb{R}^{m},
\end{equation*}
where the dot {\lq\lq}$\cdot${\rq\rq} denotes the usual inner product in $
\mathbb{R}^m$. The set of all subgradients of $f$ at ${x}$ is called
subdifferential of $f$ at ${x}$ and is denoted by $\partial f({x})$.
\end{definition}

\begin{theorem}
\label{dir} 
{\em (See \cite[Theorem 23.2]{R} or \cite[Proposition 17.7]{BC}.)}\\ 
Let $f:\mathbb{R}^m\to\mathbb{R}$ be a convex function. Then ${x}^{\ast }$ is a subgradient of $f$ at ${x}$ if and only if 
\begin{equation*}
f^{\prime }({x};{y})\geqslant {x}^{\ast }\cdot{y}%
,\quad \forall {y}\in \mathbb{R}^{m}.
\end{equation*}
\end{theorem}

\smallskip

Let us also note the following fact.

\begin{remark}
\label{W} Let 
\begin{equation*}
\Vert x\Vert _{l^{p}(\mathbb{R}^{m})}=\sqrt[p]{|x_{1}|^{p}+...+|x_{m}|^{p}}
,\qquad {x}=(x_{1},\cdots ,x_{m})\in \mathbb{R}^{m},
\end{equation*}%
denote the $l^{p}$ norm in $\mathbb{R}^{m}$. Then the ${l}^{p}$ norm is
decreasing with respect to $p\in \lbrack 1,\infty ]$. {\em This is easy to prove: given $p_{1}\leqslant p_{2}$ in $[1,\infty ]$ and ${x}=(x_{1},\cdots ,x_{m})\not=0$, let$\quad y_{i}=|x_{i}|/\Vert x\Vert _{p_{2}}$. Then $|y_{i}|\leqslant 1$, so 
\begin{equation*}
|y_{i}|^{p_{1}}\geqslant |y_{i}|^{p_{2}},\qquad \Vert {y}\Vert _{{l}
^{p_{1}}}\geqslant 1,
\end{equation*}
and, consequently, $\Vert {x}\Vert _{{l}^{p_{1}}}\geqslant \Vert {x}\Vert _{{l}^{p_{2}}}$.}
\end{remark}

\ 

The majority of the following results play a crucial role in the proof of our Theorem \ref{theorem 1}. The next theorem characterizes quite completely the subdifferential of the potential $V$ defined in \eqref{pot}-\eqref{potential}. We shall use items (a) and (c) but we find interesting by itself to furnish item (b) as well. Notice that $V$ is differentiable at matrices $X$ such that $X_\nu\not=0$ and that $X_\nu=0$ means that $X=\Omega$ if ${\nu}\not=0$ and $X_s=0$ but $X_a$ is arbitrary, if ${\nu}=0$.

\begin{theorem}
\label{convexe1} 
Let $p,q\geqslant 2$. Then, the potential $V$ defined in \eqref{pot}-\eqref{potential} is convex, differentiable at any $X\in\mathbb{R}^{d\times d}$ such that $X_\nu\not=0$, with 
\begin{equation}
\nabla V(X)=X_{\mu,p}+\widehat{\tau}_{\ast}|X_{{\nu},q}^{\widehat{ \ }}|^{-\frac{2(q-1)}{q}}X_{{\nu},q},
\label{dif V}
\end{equation}
where $X_\nu=X_s+\nu R$ (see \eqref{Bnu}), $R=X_a-\Omega$, and $X_{\mu,p}$, $X_{\nu,q}$ $X_{\nu,q}^{\widehat{ \ }}$ are defined by \eqref{B mu p}-\eqref{B hat nu q}, replacing $B$ by $X$. \ Moreover:
\medskip 

$(a)\ \quad B_{r_{\ast}}(0)\subset \partial V(\Omega)\subset
B_{\tau_\ast}(0)$, \ where $B_{r}(0)$ is the closed ball in $\mathbb{R}
^{d\times d}\equiv\mathbb{R}^{d^2}$ of radius $r$ and center $X$, and\footnote{The closed ball $B_{r_{\ast}}(0)$ is the Euclidean closed ball in $\mathbb{R}^{d\times d}$ centered at $0$ of maximum radius contained in $\partial V(\Omega)$ and $B_{\tau_\ast}(0)$ is the one of minimum radius containing $\partial V(\Omega)$.} 
\begin{equation*}
r_{\ast}=\widehat{\tau_\ast}\cdot r_q, \quad r_{q}=\left\{ 
\begin{array}{cl}
\sqrt[q]{{\nu}}/(1+{\nu}^{\frac{1}{q-2}})^{\frac{q-2}{2q}}, & \mbox{if } \ q>2\\ 
\min (1,\sqrt{{\nu}}), & \mbox{if } \ q=2;
\end{array}
\right.
\end{equation*}

$(b)\quad $ The subdifferential $\partial V(\Omega)$ has an {\lq\lq}ellipsoidal form{\rq\rq}, i.e., 
\begin{equation}
\partial V(\Omega)=\{X^*\in \mathbb{R}^{d\times d}:\ |X^*_{s}|^{q^{\prime }}+{\nu}
^{1-q^{\prime}}|X^*_a|^{q^{\prime }}\leqslant (\widehat{\tau }_{\ast})^{q^\prime}\},
\label{sub 0}
\end{equation}
$q^{\prime}=q/(q-1)$. \ {\em (See illustration below for $q=2$.)}

\medskip

$(c)\quad $ If ${\nu}=0$, the subdifferential at a matrix $\mathrm{X}$ such that $\mathrm{X}_{\nu}=0$ (i.e., $\mathrm{X}_s=0$ but $X_a$ arbitrary) is given by 
\begin{equation}
\partial V(X)=\mu_2|\mathrm{R}|^{p-2}\mathrm{R}+\{X^*\in \mathbb{R}^{d\times d}:\ X^*_a=0 \mbox{ and } |X^*_{s}|\leqslant {\tau }_{\ast }\}
\label{sub 00}
\end{equation}
(recall that $R=X_a-\Omega$).
\end{theorem}

\  

\noindent
An illustration of the subdifferential $\partial V(\Omega)$:

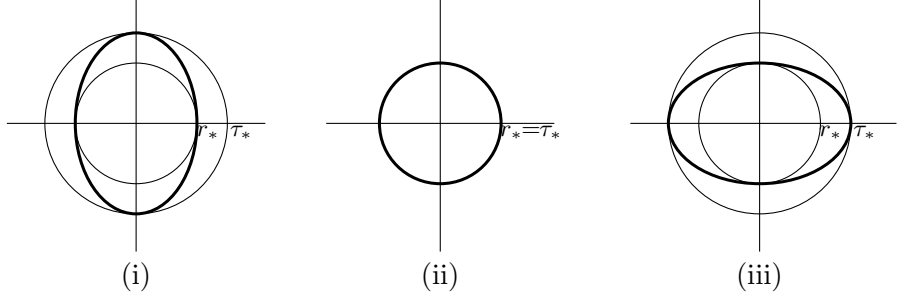
\begin{tikzpicture}[
simple/.style={draw,text width=1.5cm,align=center,minimum size=2.5cm}
]

\draw (-2.7,0) -- (0.7,0);
\draw (-1,-1.7) -- (-1,1.7);	
\draw[line width=.4mm] (-1,0) ellipse (.8cm and 1.2cm);
\draw (-1,0) circle (.8cm);
\draw (-0.34,0.1) node[below right]{\footnotesize $r_\ast$};
\draw (0.1,0.1) node[below right]{\footnotesize $\tau_\ast$};
\draw (-1,0) circle (1.2cm);
\draw (-1,-1.7) node[below]{(i)};
 
\draw (1.5,0) -- (4.5,0);
\draw (3,-1.7) -- (3,1.7);
\draw (3.65,0.1) node[below right]{\footnotesize $r_\ast\!\!=\!\!\tau_\ast$};
\draw[line width=.4mm] (3,0) circle (.8cm);
\draw (3,-1.7) node[below]{(ii)};

\draw (5.5,0) -- (9,0);
\draw (7.2,-1.7) -- (7.2,1.7);
\draw[line width=.4mm] (7.2,0) ellipse (1.2cm and .8cm);
\draw (7.2,0) circle (.8cm);
\draw (7.85,0.1) node[below right]{\footnotesize $r_\ast$};
\draw (8.3,0.1) node[below right]{\footnotesize $\tau_\ast$};
\draw (7.2,0) circle (1.2cm);
\draw (7.2,-1.7) node[below]{(iii)};
	
\end{tikzpicture}

{\small \emph{In the above picture we have the graphic in a cartesian plane $
xy$ of the the set (with the boundary in bold) given by the inequality 
$$
x^2+{\nu} ^{-1}y^2\leqslant (\widehat{\tau }_{\ast})^2,
$$
which is obtained from \eqref{sub 0} with $q=2$ by setting $x=|X^*_{s}|$ and $y=|X^*_a|$,  together with the graphics of the balls \ $x^2+y^2\leqslant
r_\ast^2$, $x^2+y^2\leqslant {\tau }_{\ast }^2$. 
We consider the three cases: (i) ${\nu}>1$; \ (ii) ${\nu}=1$ ($r_\ast=\tau_\ast$);
\ (iii) ${\nu}<1$. Notice that when $q=2$, \ \ $r_\ast=
\widehat{\tau }_{\ast }\cdot r_2=\tau_\ast\cdot\min(1,\nu)/\max(1,\nu)$.}}

\ 

\textbf{Proof} (of Theorem \ref{convexe1}){\bf{.}} \ We recall that 
$$
V(X)=U(X)+\widehat{\tau }_{\ast }W(X)
$$
with
\begin{equation}
U(X)=\frac{\mu _{1}}{p}|X_{s}|^{p}+\frac{\mu _{2}}{p}|R|^{p},\
W(X)=\sqrt[q]{|X_{s}|^q+\nu|R|^q}, \quad R=X_a-\Omega.
\label{split V}
\end{equation}

To see that $V$ is convex, we first notice that the functions 
\begin{equation*}
X\mapsto |X_{s}|^{p},\qquad X\mapsto |R|^{p}=|X_a-\Omega|^{p}
\end{equation*}
are convex since they are a composition of the convex function $t\in \mathbb{R}\mapsto |t|^{\frac{p}{2}}$ with the quadratic functions $X\mapsto
|X_{s}|^{2},$ $\ ~X\mapsto |R|^{2}$. \ In addition, we can check that
the function $W$ is convex using that any norm is convex and the fact that 
\begin{equation}
W(X)=\Vert (|X_{s}|,\sqrt[q]{{\nu}}|R|)\Vert _{{l}^{q}(\mathbb{R}^{2})}
\label{W norm}
\end{equation}
is also a norm, if ${\nu}>0$; in the case ${\nu}=0$, we have that $W(X)=|X_s|$ which is clearly convex. Thus, $V$ is convex because it is a linear combination of
convex functions.

The function $V$ is differentiable at any $X\in \mathbb{R}^{d\times d}$ such that $X_\nu\not=0$ due to the chain rule, and we can compute $\nabla V(X)$ in this case
differentiating directly $U$ and $W$ with respect to the standard variables $x_{ij}$ in $\mathbb{R}^{d\times d}$, or, alternatively, using also the chain rule one has
\begin{align}
\nabla V(X)=& \mu _{1}|X_{s}|^{p-2}X_s+\mu _{2}|R|^{p-2}R+\widehat{\tau 
}_{\ast }\nabla W(X)  \notag \\
=& X_{\mu,p}+\widehat{\tau }_{\ast }(|X_{s}|^{q}+{\nu}|R|^{q})^{\frac{1}{q}-1}(|X_{s}|^{q-2}X_{s}+{\nu} |R|^{q-2}R)  \notag \\
=& X_{\mu,p}+\widehat{\tau }_{\ast }|X_{\nu,q}^{\widehat{ \ }}|^{-\frac{2(q-1)}{q}}X_{\nu,q}.  
\label{DV}
\end{align}

Now, to show the items (a)-(c) in the statement of the theorem, let us fix a $X\in\mathbb{R}^{d\times d}$ such that $X_\nu=0$.  The function $U$ is  differentiable also at $X$, being $\nabla U(X)=0$ if ${\nu}\not=0$ and $\nabla U(X)=\mu_2|R|^{p-2}R$, if $\nu=0$. Then, the subdifferential $\partial V(X)$ is equal to $\widehat{\tau }_{\ast }\partial W(X)$ if $\nu\not=0$ and it is equal to $\mu_2|R|^{p-2}R + \widehat{\tau }_{\ast }\partial W(X)$, if $\nu=0$. Let us consider the case ${\nu}\not=0$ first -- recall that in this case, $X_\nu=0$ means that $X=\Omega$. 
The function $\underbar{W}(X)=W(X+\Omega)$ is positively homogeneous of order 1. Then by Remark~\ref{hom}, 
\begin{equation}  
W^\prime(\Omega; Y)={\underbar{W}}^{\prime }(0;Y)=\underbar{W}(Y)=W(Y+\Omega)
\label{dw}
\end{equation}
for any $Y\in \mathbb{R}^{d\times d}$. Thus, by Theorem \ref{dir}, 
\begin{equation}
X^*\in\partial W(\Omega) \ \Leftrightarrow \ W(Y+\Omega)\geqslant X^*:Y, \quad
 \forall \,Y\in \mathbb{R}
^{3\times 3}.  \label{r}
\end{equation}
Taking in the last inequality $Y=X^*$ and using Remark \ref{W}, we obtain
\begin{align*}
|X^*|^{2} \leqslant W(X^*+\Omega)=\Vert (|X^*_{s}|,\sqrt[q]{{\nu}}|X^*_{a}|)\Vert _{{l}^{q}(\mathbb{R}^{2})}\\
  \leqslant \max(1,\sqrt[q]{{\nu}})\Vert (|X^*_{s}|,|X^*_{a}|)\Vert _{{l}^{q}(\mathbb{R}^{2})} \\
 \leqslant \max (1,\sqrt[q]{{\nu}})\Vert (|X^*_{s}|,|X^*_{a}|)\Vert _{{l}^{2}(\mathbb{R}^{2})}=\max (1,\sqrt[q]{{\nu}})|X^*|,
\end{align*}
hence we have proved that 
\begin{equation}
\partial W(\mathrm{\Omega})\subset B_{\max (1,\sqrt[q]{\nu})}(0),
\label{subd W}
\end{equation}
and so the subdifferential $\partial V(\mathrm{\Omega})$ is contained in the ball $\widehat{\tau_\ast} B_{\max (1,\sqrt[q]{\nu})}(0) = B_{{\tau_\ast}}(0)$, if $\nu\not=0$.

Now, still assuming $\nu\neq0$ (and $X_\nu=0$, i.e., $X=\Omega$), let us show the claim 
\begin{equation}
B_{r_{q}}(0)\subset \partial W(\Omega)=\widehat{\tau }_{\ast
}^{-1}\partial V(\Omega)
\label{ss}
\end{equation}
of item (a), i.e., $B_{r_\ast}(0)\subset \partial V(\Omega)$ (see the definitions of $r_\ast$ and $r_q$ in item (a)). Consider an arbitrary matrix $X^*\in \mathbb{R}^{d\times
d} $ satisfying the property 
\begin{equation}
|X^*||Y|\leqslant W(Y+\Omega),\quad \forall\, Y\in \mathbb{R}^{d\times d}.  
\label{s}
\end{equation}
Then, by the Cauchy-Schwarz inequality, we have also 
\begin{equation*}
X^*:Y\leqslant W(Y+\Omega),\quad \forall \,Y\in \mathbb{R}^{d\times d}.
\end{equation*}%
By Theorem \ref{dir} we conclude that any $X^*\in \mathbb{R}^{d\times d}$ 
satisfying the property \eqref{s} belongs to $\partial W(\Omega)$. On the other
hand, by the positive homogeneity of $\underline{W}(Y)=W(Y+\Omega)$, the property \eqref{s} is equivalent to 
\begin{equation*}
|X^*|\leqslant \min_{Y}\{\underline{W}(Y):~Y\in \mathbb{R}^{d\times d} \mbox{ with }|Y|=1\}.
\end{equation*}
Let us demonstrate that this minimum is equal to $r_{q}$, which shall give
claim~\eqref{ss}. Since $1=|Y|^{2}=|Y_{s}|^{2}+|Y_{a}|^{2}$,
writing $t=|Y|_{a}^{2}$, we have $|Y_{s}|^{2}=1-t$ and
\begin{equation*}
\underline{W}(Y)\equiv \underline{W}(t)=\sqrt[q]{(1-t)^{\frac{q}{2}}+{\nu}t^{\frac{q}{2}}},\quad
t\in \lbrack 0,1].
\end{equation*}%
By a straightforward computation, we obtain that the minimum of the function 
$\alpha (t)=(1-t)^{\frac{q}{2}}+{\nu}t^{\frac{q}{2}}$ in the interval $[0,1]$
is $(r_{q})^{q}$
. Thus, we have proved claim (a).

To show claim (b), we follow a similar reasoning as above, with the help of the H\"{o}lder inequality. From \eqref{dw} and Theorem \ref{dir} we have that 
\begin{equation}  \label{xsub}
X^*\in \partial W(\Omega) \ \Leftrightarrow \ X^*:Y\leqslant W(Y+\Omega),\ \forall
\,Y\in \mathbb{R}^{d\times d}.
\end{equation}
Taking $Y$ in this inequality with $Y_{s}=|X^*_{s}|^{q^{\prime }-2}X^*_{s}$ and $Y_{a}={\nu}^{1-q^{\prime}}|X^*_{a}|^{q^{\prime }-2}X^*_{a}$, we obtain 
\begin{equation*}
|X^*_{s}|^{q^{\prime }}+{\nu}^{1-q^{\prime}}|X^*_{a}|^{q^{\prime
}}\leqslant \sqrt[q]{|X^*_{s}|^{q^{\prime }}+{\nu}^{1-q^{\prime}}|X^*_{a}|^{q^{\prime}}},
\end{equation*} 
hence $|X^*_{s}|^{q^{\prime }}+{\nu}^{1-q^{\prime}}|X^*_{a}|^{q^{\prime}}\leqslant 1$. Reciprocally, if this inequality holds then, using the H\"{o}lder inequality, we have 
\begin{eqnarray*}
X^* &:&Y=X^*_{s}:Y_{s}+({\nu}^{-1/q}X^*_{a}):(\nu^{1/q}Y_{a}) \\
&\leqslant &\sqrt[q^{\prime }]{|X^*_{s}|^{q^{\prime }}+{\nu}^{-q^{\prime}/q}|X^*_{a}|^{q^\prime}}\cdot\sqrt[q]{|Y_{s}|^{q}+{\nu}|Y_{a}|^{q}}\leqslant W(Y+\Omega),
\end{eqnarray*}
for any $Y\in \mathbb{R}^{d\times d}$. Therefore, again by Theorem \ref{dir}, we obtain that $X^*\in \partial W(\Omega)$.

In the case ${\nu}=0$, as we already observed above, we have that $W(Y)=~|Y_s|$. Thus, by the definition \ref{subd} of subdifferential and $X_\nu=0$ (i.e., $X_s=0$),
$$
X^*\in\partial W(X) \ \Leftrightarrow \ |Y_s| \ge X^*:(Y-X_a), \quad \forall\, Y\in\mathbb{R}^{d\times d}.
$$
Taking $Y=X^*_a+X_a, \ X^*_s+X_a$, we obtain that $X^*_a=0$ and $|X^*_s|\le 1$. Reciprocally, if $X^*_a=0$ and $|X^*_s|\le 1$ then by Cauchy-Schwarz
inequality we have that $X^*:(Y-X_a)=X^*_s:Y_s\le |Y_s|$ for all $Y\in\mathbb{R}
^{d\times d}$, thus $X^*\in\partial W(X)$. \ \ $\blacksquare$

\bigskip

As a corollary of the above computations we have the following result that justifies the definition of $\widehat{\tau}_{\ast}$, which was given in \eqref{tau star hat}.

\begin{corollary}
\label{est DW}
For any matrix $X\in \mathbb{R}^{d\times d}$ such that $X_\nu\not=0$
one has the estimate 
\begin{equation}
\frac{|X_{{\nu},q}|}{|X_{{\nu},q}^{\widehat{ \ }}|^{\frac{2(q-1)}{q}}}\leqslant \max
(1,\sqrt[q]{\nu}).  
\label{1c}
\end{equation}
\end{corollary}

\textbf{Proof}. \ As in the computations in \eqref{DV}, we have that 
\begin{equation*}
|X_{{\nu},q}^{\widehat{ \ }}|^{-\frac{2(q-1)}{q}}{X_{{\nu},q}}=\nabla W(X),
\end{equation*}
recalling that $W$ was defined in \eqref{split V}. Now we claim that $\nabla W(X)\in \partial W(\Omega)$, and thus the estimate \eqref{1c} shall follow from this fact, since $\partial W({\Omega})\subset B_{\max (1,\sqrt[q]{\nu})}(0)$ -- see \eqref{subd W}. Accounting \eqref{dw}-\eqref{r}, it remains then only to show that $\nabla W(X):Y\leqslant \underline{W}(Y):=W(Y+\Omega)$ for all $Y\in \mathbb{R}^{d\times d}$. But as $W$ is differential at $X$, $\nabla W(X):Y$ coincides with the one-sided directional derivative $W'(X;Y)=\underline{W}(X-\Omega; Y)$, which by Theorem \ref{dir der} is given by
\begin{equation*}
\inf_{\lambda >0}\lambda ^{-1}\left( \underline{W}(X-\Omega+\lambda Y)-\underline{W}(X-\Omega)\right) .
\end{equation*}%
In addition, $\underline{W}$ is a seminorm, then 
\begin{equation*}
\underline{W}(X-\Omega+\lambda Y)\leqslant \underline{W}(X-\Omega)+\lambda\underline{W}(Y)
\end{equation*}%
Summing up, we obtain our above claim.$\hfill \;\blacksquare $

\medskip

The next corollary says that the condition $\mathbb{S}\in\partial V$ implies that the relation \eqref{SR01g} is satisfied, if $\nu>0$ or $\nu=\mu_2=0$.

\begin{corollary}
\label{S}
Let $\mathbb{S}$ be a tensor in $\mathbb{R}^{d\times d}$ such that 
$\mathbb{S}\in \partial V(\mathrm{B})$ (more precisely, $\mathbb{%
S}(x,t)\in \partial V_{(x,t)}(\mathrm{B}(x,t))$, i.e., satisfies the variational inequality 
\begin{equation*}
V_{(x,t)}(X)-V_{(x,t)}(\mathrm{B}(x,t))\geqslant \mathbb{S}(x,t):(X-\mathrm{B%
}(x,t)),\quad \forall \ X\in \mathbb{R}^{d\times d},
\end{equation*}
for almost all $(x,t)\in\mathcal{O}_T$). Then $\mathbb{S}$ accomplishes the relation \eqref{SR01g} for almost all $(x,t)\in\mathcal{O}_T$, when $\nu>0$ or $\nu=\mu_2=0$.
\end{corollary}

\textbf{Proof}. \ If $\mathrm{B}_\nu\not=0$ then, by the first claim in Theorem \ref{convexe1}, we have that 
\begin{equation*}
\partial V(\mathrm{B})=\nabla V(\mathrm{B})=\mathrm{B}_{\mu,p}+\widehat{\tau}_{\ast }{|\mathrm{B}_{{\nu},q}^{\widehat{ \ }}|^{-\frac{2(q-1)}{q}}}{\mathrm{B}_{{\nu},q}}.
\end{equation*}
Then $\mathbb{S}=\nabla V(\mathrm{B})$ and satisfies \eqref{SR01g}. On the other hand, if $\mathrm{B}_\nu=0$ or $\nu=\mu_2=0$ then $|\mathbb{S}|\le \tau_{\ast}$, by items (a) and (c) in Theorem \ref{convexe1}. $\hfill \;\blacksquare $

\medskip

\begin{corollary}
\label{mono} The modified plastic operator $\mathbf{\mathcal{P}}_{{\nu},q}$, introduced in \eqref{mpo}, is monotone, for any values of ${\nu}>0$ or $\nu=\mu_2=0$ and $q\ge2$.
\end{corollary}

\textbf{Proof}. \ By the proof of Theorem \ref{convexe1} we have that $\mathbf{\mathcal{P}}_{{\nu},q}(\mathrm{X}_{s},\mathrm{X}_{a})\in \partial W(\mathrm{X})$,  where $W$ is a convex function. This ends the proof because (as it is easy to check) any subgradient of a convex function is a monotone operator. $\hfill
\;\blacksquare $

\bigskip

Next we present some technical results we shall use latter.

\begin{lemma}
Let $W=W(X)$ be a positive convex function on $X\in \mathbb{R}^{d\times d}$.
\ Then, for any $q\geqslant 2$ \ and natural $n$, the approximated function 
\begin{equation*}
W_{n}(X)=\sqrt[q]{\left( W(X)\right) ^{q}+n^{-1}}
\end{equation*}%
is also convex with respect to the parameter $X\in \mathbb{R}^{d\times d}$.
\end{lemma}

\textbf{Proof}. \ Note that the function $\varphi (z)=\sqrt[q]{z^{q}+n^{-1}%
}$ is monotone increasing and convex function for $z\geqslant 0.$ \
Therefore \ applying the definition of convex function, we easily derive
that the composition $W_{n}(X)=\varphi (W(X))$ is also convex with respect
to the parameter $X\in \mathbb{R}^{d\times d}$. $\hfill \;\blacksquare $

\bigskip

An auxiliary result that we will often use is the equivalence of the norms 
\begin{equation*}
\Vert \mathbf{v}\Vert _{V_{p}}=\Vert \mathbf{v}\Vert _{L^{2}(\mathcal{O} )}+\Vert
\mathrm{B}\Vert _{L_{p}(\mathcal{O})},\quad \Vert \mathbf{v}\Vert
_{W^{1,p}(\mathcal{O})}=\Vert\mathbf{v}\Vert _{L^{p}(\mathcal{O} )}+\Vert\mathrm{B}\Vert _{L^{p}(\mathcal{O})},
\end{equation*}
for $p\ge2$, where $W^{1,p}(\mathcal{O})$ denotes the usual Sobolev space (recall that we have set $\mathrm{B}=\nabla\mathbf{v}$).
Since $p\ge2$, we have that $\Vert\mathbf{v}\Vert _{L^{2}(\mathcal{O})}\le |\mathcal{O}|^{(p-2)/2p}\Vert\mathbf{v}\Vert _{L^{p}(\mathcal{O})}$ (where $|\mathcal{O}|$ is the Lebesgue measure of $\mathcal{O}$). Thus, the above equivalence is a direct   consequence of the following lemma.

\begin{lemma}
\label{equivalents} 
For any $p\geqslant 2$ there exists a constant $C$ such that
\begin{equation}
\Vert \mathbf{v}\Vert_{L^p(\mathcal{O})}\leqslant C\Vert \mathbf{v}\Vert
_{V_{p}},\quad \forall \ \mathbf{v}\in W^{1,p}(\mathcal{O}).
\label{eq}
\end{equation}
\end{lemma}

\textbf{Proof.} \ First we observe the continuous embeddings $W^{1,p}(\mathcal{O})\subset L^{p}(\mathcal{O})\subset L^{2}(\mathcal{O})$ being compact the first one. Then by Ehrling's lemma (also known as Lions' lemma - see \cite[Exercise 6.12, p. 173]{Brezis}) for any $\epsilon>0$ there exists a constant $C=C(\epsilon)$ such that 
\begin{equation}
\Vert \mathbf{v}\Vert _{L^p(\mathcal{O})}\leqslant \epsilon\Vert\mathbf{v}\Vert
_{W^{1,p}(\mathcal{O})}+C\Vert\mathbf{v}\Vert_{L^2(\mathcal{O})},\quad \forall \ \mathbf{v}\in W^{1,p}(\mathcal{O}).
\label{Ehrling}
\end{equation}
Choosing $\epsilon=1/2$ in this inequality we obtain the desired result. $\hfill \;\blacksquare$

\bigskip

In the same vein we can slightly improve the Korn's inequality
\begin{equation}
\Vert\mathrm{B}\Vert _{L^p(\mathcal{O})}\leqslant C\left(\Vert \mathbf{v}\Vert
_{L^p(\mathcal{O})} + \Vert\mathrm{B}_s\Vert _{L^p(\mathcal{O})}\right), \quad p\in (1,\infty), \ C=C(\Omega)
\label{Korn p}
\end{equation}
(see \cite[(1.70), p. 20]{tem2}) to obtain

\begin{proposition} 
\label{Korn}
For any $p\in [2,\infty)$ there exists a constant $C$ such that
$$
\Vert\mathrm{B}\Vert _{L^p(\mathcal{O})}\leqslant C\left(\Vert \mathbf{v}\Vert
_{L^2(\mathcal{O})} + \Vert\mathrm{B}_s\Vert _{L^p(\mathcal{O})}\right),\quad \forall \ \mathbf{v}\in W^{1,p}(\mathcal{O}).
$$
\end{proposition}

\textbf{Proof.}\footnote{With some restriction on $p$, a different proof can be made by using the interpolation inequality $\Vert\mathbf{v}\Vert_{L^p(\mathcal{O})}\le C\Vert\mathbf{v}\Vert_{L^2(\mathcal{O})}^a\Vert\mathrm{B}\Vert_{L^2(\mathcal{O})}^{1-a}$, $a\in [0,1]$, in \eqref{Korn p}, and then the {\lq}Young inequality with $\epsilon${\rq}.} \ For $\epsilon<1$, the inequality \eqref{Ehrling} can be rewritten as
$$
\Vert \mathbf{v}\Vert _{L^p(\mathcal{O})}\leqslant \frac{\epsilon}{1-\epsilon}\Vert\mathrm{B}\Vert_{L^p(\mathcal{O})}+\frac{C}{1-\epsilon}\Vert\mathbf{v}\Vert_{L^2(\mathcal{O})}.
$$
Substituing this estimation in \eqref{Korn p}, we obtain
$$
\Vert\mathrm{B}\Vert_{L^p(\mathcal{O})}\leqslant C\left(\frac{\epsilon}{1-\epsilon}\Vert\mathrm{B}\Vert_{L^p(\mathcal{O})} + \frac{C}{1-\epsilon}\Vert\mathbf{v}\Vert_{L^2(\mathcal{O})} + \Vert\mathrm{B}_s\Vert_{L^p(\mathcal{O})}\right) 
$$
with possibly different $C$'s here. Hence, choosing $\epsilon$ such that $C\epsilon/(1-\epsilon)=1/2$ we end the proof. \ $\blacksquare$

\bigskip

Let us remember the well-known interpolation result \cite{Evans}.

\begin{lemma}
\label{interpolation} Let $s_{1},s_{2},r\in \lbrack 1,+\infty ]$ \ and $\
\theta _{1},\theta _{2}\in \lbrack 0,1]$ such that\break\hfill \ $\theta _{1}+\theta _{2}=1,\quad \frac{\theta _{1}}{s_{1}}+\frac{\theta _{2}}{s_{2}}=\frac{1}{r}$. \ Then 
\begin{equation}
\Vert \mathbf{v}\Vert_{L^{r}(\mathcal{O})}\leqslant\Vert\mathbf{v}\Vert_{L^{s_{1}}(\mathcal{O})}^{\theta_{1}}\Vert\mathbf{v}\Vert_{L^{s_{2}}(\mathcal{O})}^{\theta_{2}},\qquad \forall \mathbf{v}\in L^{s_{1}}(\mathcal{O})\cap L^{s_{2}}(\mathcal{O}). \label{interp}
\end{equation}
\end{lemma}

Since $\mathcal{O}$ is a bounded domain with $C^{2}$-smooth boundary, the
Gagliardo--Nirenberg--Sobolev embedding theorem (see \cite{Evans}) and Lemma 
\ref{equivalents} imply the following result.

\begin{lemma}
\label{embedding} 
There exists a positive constant $C$, such that
\begin{equation}
\Vert\mathbf{v}\Vert_{L^{r}(\mathcal{O})}\leqslant C\Vert \mathbf{v}\Vert
_{V_{p}},\quad \forall \mathbf{v}\in V_{p},  \label{embed}
\end{equation}
with $r=dp/(d-p)$, if $1\le p<d$; any $r<\infty$, if $p=d$; and $r=\infty$, if $p>d$.
\end{lemma}

\section{Construction of the approximation problem}
\label{Existence}

In this section we consider an approximated problem for the system \eqref{OneDeq1}. We
solve this approximated problem applying the Faedo-Galerkin method (see for instance \cite{AKM}).

\medskip

Since the space $V_{p}$ is separable, it is the span of a countable set of
linearly independent functions $\left\{ \mathbf{e}_{k}\right\}
_{k=1}^{\infty }$. 
For simplicity of calculations, we choose this set as the eigenfunctions for
the problem 
\begin{equation}
(\mathbf{e}_{k},{\boldsymbol{\varphi }})_{W}=\lambda _{k}(\mathbf{e}_{k},{%
\boldsymbol{\varphi }}),\qquad \forall {\boldsymbol{\varphi }}\in
W=W_{2}^{3}(\mathcal{O} )\cap V,  \label{wh}
\end{equation}%
where $(\cdot ,\cdot )$ and $(\cdot ,\cdot )_{W}=(\cdot ,\cdot
)_{W_{2}^{3}(\mathcal{O} )}$ are the standard inner products in $L^{2}(\mathcal{O} )$
and $W_{2}^{3}(\mathcal{O} ),$ respectively. The solvability of the problem %
\eqref{wh} follows from the spectral theory (see \cite{Boyer}, \cite{Evans}%
).\ \ This theory permits to construct this set $\left\{ \mathbf{e}%
_{k}\right\} _{k=1}^{\infty }$ as an orthogonal basis for $W$ and an
orthonormal basis for $H$.

We can consider the subspace $W^{n}=\mathrm{span}\,\{\mathbf{e}_{1},\ldots ,%
\mathbf{e}_{n}\}$ of $W$, for any fixed natural $n.$ Let $P_{n}:W\rightarrow
W_{n}$ be the orthogonal projection defined by%
\begin{equation*}
P_{n}{\boldsymbol{\varphi }}=\sum_{k=1}^{n}\widetilde{\beta }_{k}\widetilde{%
\mathbf{e}}_{k}=\sum_{k=1}^{n}\beta _{k}\mathbf{e}_{k},\qquad \forall {%
\boldsymbol{\varphi }}\in W,
\end{equation*}%
with $\widetilde{\beta }_{k}=\left( {\boldsymbol{\varphi }},\widetilde{%
\mathbf{e}}_{k}\right) _{W}\ $and $\beta _{k}=\left( {\boldsymbol{\varphi }},%
\mathbf{e}_{k}\right) $ and $\{\widetilde{\mathbf{e}}_{j}=\frac{1}{\sqrt{%
\lambda _{j}}}\mathbf{e}_{j}\}_{j=1}^{\infty }$ is the orthonormal basis of $%
\ W.$ \ By Parseval's identity, we get 
\begin{eqnarray}
||P_{n}{\boldsymbol{\varphi }}||_{H} &\leqslant &||{\boldsymbol{\varphi }}%
||_{H},\qquad ||P_{n}{\boldsymbol{\varphi }}||_{W}\leqslant ||{\boldsymbol{%
\varphi }}||_{W},  \notag \\
P_{n}{\boldsymbol{\varphi }} &\longrightarrow &{\boldsymbol{\varphi }}\quad 
\mbox{ strongly in
	}\text{ }H \, \mbox{ and strongly in
	} \text{ }W,  \label{2}
\end{eqnarray}%
where 
\begin{equation*}
\Vert {\boldsymbol{\varphi }}\Vert _{H}=\sqrt{({\boldsymbol{\varphi }},{%
\boldsymbol{\varphi }})}\quad \text{and}\quad \Vert {\boldsymbol{\varphi }}%
\Vert _{W}=\sqrt{({\boldsymbol{\varphi }},{\boldsymbol{\varphi }})_{W}}
\end{equation*}%
are the norms in $H$ and $W$, respectively. Obviously, we have $W\subset
V_{p}\subset H$.

Considering an arbitrary ${\boldsymbol{\varphi }}\in L^{p}(0,T;W),$ we also
have%
\begin{equation*}
||P_{n}{\boldsymbol{\varphi }}(t)||_{W}\leqslant ||{\boldsymbol{\varphi }}%
(t)||_{W}\quad \mbox{ and
	}\quad P_{n}{\boldsymbol{\varphi }}(t)\rightarrow {\boldsymbol{\varphi }}%
(t)\quad \text{ strongly in}\ W,
\end{equation*}%
which are valid for a.e. $t\in (0,T).$ Hence, Lebesgue's dominated
convergence theorem implies that for any ${\boldsymbol{\varphi }}\in
L^{p}(0,T;W)$, we obtain 
\begin{equation}
P_{n}{\boldsymbol{\varphi }}\longrightarrow {\boldsymbol{\varphi }}\quad 
\mbox{ strongly in
	}\ L^{p}(0,T;W)\quad \text{as }n\rightarrow \infty .  \label{c02}
\end{equation}

\medskip

Now let us define the vector function
\begin{equation}
\mathbf{v}^{n}(t)=\sum_{k=1}^{n}c_{k}^{(n)}(t)\ \mathbf{e}_{k},\qquad
c_{k}^{(n)}(t)\in {\mathbb{R}},  \label{vv}
\end{equation}
as the solution of the system
\begin{equation}
\left\{ 
\begin{array}{l}
\int\limits_{\mathcal{O} }[\partial _{t}\mathbf{v}^{n}\mathbf{e}_{k}+\left( 
\mathbf{v}^{n}\mathbf{\cdot \nabla }\right) \mathbf{v}^{n}\mathbf{e}_{k}+%
\mathbb{S}^{n}:\nabla \mathbf{e}_{k}]\,d{x} \\ 
\qquad \qquad +\int_{\Gamma}(\nabla g(\mathbf{v}^{n})\cdot{\bm{\tau}})(\mathbf{e}_{k}\cdot{\bm{\tau}})\,d{\Gamma }\,dt=0,\qquad \forall
k=1,2,\dots ,n, \\ 
\mathbf{v}^{n}(0)=\mathbf{v}_{0}^{n}.%
\end{array}%
\right.  \label{y1}
\end{equation}
The matrix functions $\mathbb{T}^{n}$, $\mathbb{S}^{n}$ are prescribed by the relations
\begin{equation}
\mathbb{T}^n=-p^n\,\mathrm{I}+\mathbb{S}^n,\qquad\mathbb{S}^n=\mathrm{B}_{\mu,p}^n+\widehat{\tau}_{\ast}\frac{\mathrm{B}_{\nu,q}^n}{\sqrt[q]{(|{\mathrm{B}_{\nu,q}^{n\,{\widehat{ \ }}}|^2+n^{-1})^{q-1}}}} \quad \mbox{(cf. \eqref{approximate stress})}
\label{delta}
\end{equation}
and the matrix functions $\mathrm{B}^n=\nabla\mathbf{v}^{n}$, $\mathrm{B}_{\mu,p}^{n},\,\mathrm{B}_{\nu,q}^{n},\,\mathrm{B}_{\nu,q}^{n}$ are calculated through the formulas \eqref{B mu p}, \eqref{B nu q}, \eqref{B hat nu q}. The function $\mathbf{v}_{0}^{n}$ \ is the orthogonal projection of $\mathbf{v}_{0}\in H$\ into the space $W^{n}$. Note that the system \eqref{vv}-\eqref{delta} is a weak formulation of the problem \eqref{approximate ibvp}. 

Next we prove that the approximated problem \eqref{y1} is solvable. As the starting point let us prove Lemma \ref{apr copy(1)}.

\smallskip

\textbf{Proof of Lemma \ref{apr copy(1)}.} \ Item (b) is obtained by straightforward computation, the definition \eqref{approximate stress} of $\mathbb{S}^n$ and Corollary \ref{est DW}.  Regarding item (a), in a matrix $X$ such that $X_\nu\not=0$, $V$ is differentiable and $V'(X;X)=\nabla V(X):X$. Then by Theorem \ref{convexe1}/\eqref{dif V} we have
\begin{equation}
\begin{array}{rl}
V'(X;X)&=X_{\mu,p}:X+\widehat{\tau}_{\ast}|X_{{\nu},q}^{\widehat{ \ }}|^{-\frac{2(q-1)}{q}}X_{{\nu},q}:X\\
&=\mu_1|X_s|^p+\mu_2|R|^{p-2}(R:X_a)\\
&\qquad\qquad + \ \ \widehat{\tau}_{\ast}|X_{{\nu},q}^{\widehat{ \ }}|^{-\frac{2(q-1)}{q}}\left(|X_s|^q+\nu |R|^{q-2}(R:X_a)\right).
\end{array}
\label{est V'}
\end{equation}
Now, $R:X_a=(X_a-\Omega):X_a=|X_a|^2-(\Omega:X_a)\ge |X_a|^2 - |\Omega|\,|X_a|=|X_a|\left(|X_a| - |\Omega|\right)$, then, if $|X_a|\ge |\Omega|$, from \eqref{est V'} we obtain that $V'(X;X)\ge \mu_1|X_s|^p$. On the contrary, i.e., if $|X_a|<|\Omega|$, we use the following estimates in \eqref{est V'}:
$$
|R|^{p-2}|(\Omega:X_a)|\le 2^{p-3}(|X_a|^{p-2}+|\Omega|^{p-2})|\Omega|\,|X_a|\le 2^{p-2}|\Omega|^p;
$$
$$
\begin{array}{rl}
&\widehat{\tau}_{\ast}|X_{{\nu},q}^{\widehat{ \ }}|^{-\frac{2(q-1)}{q}}\nu |R|^{q-2}|(R:X_a)| = \widehat{\tau}_{\ast}{\displaystyle\frac{\nu |R|^{q-2}|(R:X_a)|}{(|X_s|^q+\nu|R|^q)^{\frac{q-1}{q}}}}\\
&\le \widehat{\tau}_{\ast}{\displaystyle\frac{\nu |R|^{q-2}|(R:X_a)|}{(\nu|R|^q)^{\frac{q-1}{q}}}} \le \widehat{\tau}_{\ast}\sqrt[q]{\nu} |X_a| \le \tau_\ast |\Omega|.
\end{array}
$$
Then, disregarding non negatives terms in \eqref{est V'}, we arrive at the estimate from below in \eqref{estimate for V and Vn} for $V$, at those $X$ such that $X_\nu\not=0$. Next, estimating \eqref{est V'} from above, using Corollary \ref{est DW}, we have
$$
\begin{array}{rl}
V'(X;X)&\le \mu_1|X|^p+\mu_2 |R|^{p-1}|X| + \tau_\ast |X| \\
       &\le \mu_1|X|^p+\mu_2 2^{p-2}(|X|^{p-1}+|\Omega|^{p-1})|X| + \tau_\ast |X|\\
			 &\le \mu_1|X|^p+\mu_2 2^{p-2}\left(|X|^p+\frac{1}{p}|X|^p + (1-\frac{1}{p})|\Omega|^p\right) + \tau_\ast |X|,
\end{array} 
$$
where for the last inequality we have used the Young's inequality $ab\le\frac{a^p}{p}+\frac{b^{p'}}{p'}$, \ {\small $p'=p/(p-1)$}. Thus, we arrived at the estimate \eqref{estimate for V and Vn} from above for $V$, at any $X$ such that $X_\nu\not=0$.

In the case of $X_\nu=0$, we observe that the potential $U$ is differentiable at any matrix $X\in\mathbb{R}^{d\times d}$ and the above computations shows the estimate
\begin{equation}
\mu_1 |X_s|^p - \mu_2 2^{p-2}|\Omega|^p \ \le \ V'(X;X)  \ \le \ c_1|X|^p + c_2|\Omega|^p, \qquad \forall\ X\in \mathbb{R}^{d\times d}.
\label{est U'}
\end{equation}
Then, to end estimate \eqref{estimate for V and Vn} for $V$ it remains only to estimate $\widehat{\tau_\ast}W'(X,X)$ in the case of $X_\nu=0$. But, in this case it is easy to check, using the definitions of $X_\nu$ (i.e., $X_\nu=X_s+\nu R$, $R=X_a-\Omega$) and $W$, and the definition, given in \eqref{def f'}, of one-sided directional derivative, that $W(X;X)=0$ for all $X$ such that $X_\nu=0$. Thus, we have ended the proof of \eqref{estimate for V and Vn} for the potential $V$.

With respect to $V^n$, by the fact that the only difference between $V$ and $V^n$ is in $W$ modified to $W^n$, to obtain \eqref{estimate for V and Vn} for the potential $V^n$ we need only to estimate $\widehat{\tau_\ast} W^n(X;X)$. This estimation is even easier than the estimate we have made to $W$, for $W^n$ is differentiable at any $X\in\mathbb{R}^{d\times d}$ and it is very similar to $W$. Indeed, for every $X\in\mathbb{R}^{d\times d}$, we have
$$
\begin{array}{rl}
W^n(X;X)&=\nabla W^n(X):X={\displaystyle\frac{\mathrm{X}_{{\nu},q}}{\sqrt[q]{(|\mathrm{X}_{{\nu},q}^{\widehat{ \ }}|^{2}+n^{-1})^{q-1}}}}:X\\
        &={\displaystyle\frac{|\mathrm{X}_s|^q + \nu |R|^{q-2}(R:X_a)}{\left(|\mathrm{X}_s|^q + \nu |\mathrm{R}|^q+n^{-1}\right)^{\frac{q-1}{q}}}},
\end{array}
$$
then if follows the estimate from above
$$
W^n(X;X)\le {\displaystyle\frac{\left||\mathrm{X}_s|^q + \nu |R|^{q-2}(R:X_a)\right|}{\left(|\mathrm{X}_s|^q + \nu |\mathrm{R}|^q\right)^{\frac{q-1}{q}}}}
$$
and then we can proceed exactly as above to get the same estimate we have obtained for $\widehat{\tau_\ast}W^(X;X)$, i.e.,
$$
\widehat{\tau_\ast}W^n(X;X)\le \tau_\ast |X|.
$$
The estimation of $\widehat{\tau_\ast}W^n(X;X)$ from below does not reduce to the previous estimation we have made for $\widehat{\tau_\ast}W(X;X)$ but we can proceed in a similar way, first observing as above that if $|X_a|\ge |\Omega|$ then $W^n(X;X)\ge 0$, so we can completely disregard $W^n(X,X)$ in this case, and if $|X_a|<|\Omega$, disregarding non negative terms we have, similarly as above:
$$
\begin{array}{rl}
\widehat{\tau_\ast}W^n(X;X)&\!=\,\widehat{\tau_\ast}{\displaystyle\frac{|\mathrm{X}_s|^q + \nu |R|^{q-2}(R:X_a)}{\left(|\mathrm{X}_s|^q + \nu |\mathrm{R}|^q+n^{-1}\right)^{\frac{q-1}{q}}}}\\
&\!\ge\, \widehat{\tau_\ast}{\displaystyle\frac{\nu |R|^{q-2}(R:X_a)}{\left(|\mathrm{X}_s|^q + \nu |\mathrm{R}|^q+n^{-1}\right)^{\frac{q-1}{q}}}}\\
&\!\ge\, - \widehat{\tau_\ast}\sqrt[q]{\nu}|X_a| \ge - \tau_\ast |\Omega|.
\end{array}
$$
$\hfill \;\blacksquare$

\bigskip\bigskip

\begin{lemma}
\label{1} Let the conditions on $\Omega$ and $g$ stated in Theorem \ref{theorem 1} be in force. Then for any $\mathbf{v}_0\in H$, $p\geqslant 2$ and $q\geqslant 2$, there exists a solution 
\begin{equation*}
\mathbf{v}^{n}\in L^{\infty}(0,T;H)\cap L^{p}(0,T;V_{p})
\end{equation*}
of the system \eqref{vv}-\eqref{delta} that satisfies the following estimates:
\begin{equation}
\int\limits_{\mathcal{O}}|\mathbf{v}^{n}(x,t)|^{2}d{x}+\int\limits_{\mathcal{O}_t}|{\nabla }\mathbf{v}^{n}|^{p}d{x}\,dt
\leqslant C,\quad t\in \lbrack 0,T],  
\label{omegadelta}
\end{equation}
\begin{equation}
\Vert \mathbb{S}^{n}\Vert _{L^{p^{\prime}}(\mathcal{O}_{T})}\leqslant C
\label{Sdelta}
\end{equation}
and
\begin{equation}
\Vert\partial_{t}\mathbf{v}^{n}\Vert_{L^{p^{\prime}}(0,T;\,W^{\ast})}\leqslant C.  \label{timedelta}
\end{equation}
The positive constants $C$ \textit{do not depend on} $n$, but may depend on $p$, $q$, $\mathbf{v}_{0}$, $\mu _{i}$, $\nu$ and the boundary function $g$. \ (Recall that $\mathcal{O}_t=\mathcal{O}\times (0,t)$.)
\end{lemma}

\textbf{Proof}.\ \ We have that $c_{k}^{(n)}(t)=\int\limits_{\mathcal{O} }\mathbf{%
ve}_{k}d{x}$ for $k=1,...,n$ with $c_{k}^{(n)}$ introduced in \eqref{vv}. Therefore the system \eqref{y1} is a system of $n$ nonlinear ordinary differential equations of the first order, which can be written in
the form
\begin{eqnarray*}
\frac{d\mathbf{c}^{(n)}}{dt} &=&F(\mathbf{c}^{(n)},t),\qquad t\in \left[ 0,T%
\right] , \\
\mathbf{c}^{(n)}(0) &=&(c_{1}^{(n)}(0),...,c_{n}^{(n)}(0)),\qquad
c_{k}^{(n)}(0)=\int\limits_{\mathcal{O} }\mathbf{v}_{0}\mathbf{e}_{k}d\mathbf{x,}
\end{eqnarray*}
for the vector function $\mathbf{c}^{(n)}(t)=(c_{1}^{(n)}(t),...,c_{n}^{(n)}(t))$. The function $F(\mathbf{c}^{(n)},t)$ is a continuous function in $\mathbf{c}^{(n)}$ for earch fixed $t$ and $F(\mathbf{c}^{(n)},t)$ is measurable in $t$ for each fixed $\mathbf{c}^{(n)}$. Therefore Carath\'{e}odory's existence theorem (a generalization of Peano's theorem) of the theory of ordinary differential equations implies the existence of a local time solution $\mathbf{c} ^{(n)}\in C^{1}([0,t_{n}))$ for a maximal existence time interval $[0,t_{n}), $ such that $t_{n}\leqslant T$.

Let us demonstrate that $t_{n}=T.$ To do it we deduce a priori estimates for 
$\mathbf{c}^{(n)}.$ Let us multiply \eqref{y1}$_{1}$ by $c_{k}^{(n)}$ and
take the sum on the index $k=1,...,n.$\ Then the integration over the time
interval $(0,t)$ gives the {\lq\lq}energy equation{\rq\rq} (cf. \eqref{energy eq})
\begin{equation}
\begin{array}{c}
{\textstyle\frac{1}{2}}{\displaystyle\int_{\mathcal{O}}|\mathbf{v}^{n}(x,t)|^{2}d{x}+\int_{\mathcal{O}_t}\mathbb{S}^{n}:{\nabla}\mathbf{v}^{n}d{x}\,dt+\int_{\Gamma_t}\nabla (g(\mathbf{v}^{n})\cdot\bm{\tau})(\mathbf{v}^{n}\cdot\bm{\tau})\,d{\Gamma}\,dt}\\
={\textstyle\frac{1}{2}}{\displaystyle\int_{\mathcal{O} }|\mathbf{v}_{0}^{n}|^{2}d{x}}.
\end{array}
\label{za}
\end{equation}
The estimate \eqref{estimate for V and Vn} for $V^n$ and the condition $\nabla g(\mathbf{v})\cdot\mathbf{v}\ge 0$ for every $\mathbf{v}\in\mathbb{R}^3$ (see \eqref{condition for g}) imply from here 
\begin{equation*}
{\textstyle\frac{1}{2}}\int_{\mathcal{O}}|\mathbf{v}^{n}|^{2}d{x}+\mu_1\int_{0}^{t}\int_{\mathcal{O}}|\mathrm{B}^{n}_s|^pd{x}\,dt
\le C + {\textstyle\frac{1}{2}}\int_{\mathcal{O} }|\mathbf{v}_{0}^{n}|^{2}d{x},
\end{equation*}
where $C=\int_{\mathcal{O}_T}[\,\mu_2 2^{p-2}|\Omega|^p - \tau_\ast |\Omega|\,]\,d{x}\,dt$. Then using the Korn's inequality/Proposition \ref{Korn}, we deduce the a priori estimate \eqref{omegadelta}.

Now we are able to show that $t_{n}=T.$ If $t_{n}<T,$ \ since $[0,t_{n})$ is
the maximal time interval, then $|\mathbf{c}^{(n)}(t)|^{2}=\Vert \mathbf{v}%
^{n}(t)\Vert _{L^{2}(\mathcal{O} )}^{2}$ has to tend to $+\infty $ as $%
t\rightarrow t_{n}$ which contradicts \eqref{omegadelta}.

From the estimate \eqref{est Sn} we obtain the a priori estimate \eqref{Sdelta}. 

Let ${\boldsymbol{\varphi }}\in C_{0}^{1}(0,T;W)$ and $P_{n}{\boldsymbol{%
\varphi }}$ be an arbitrary function and its orthogonal projection of $W$\
onto $W^{n}.$\ The first equality of \eqref{y1} is linear with respect of
the functions $\mathbf{e}_{k}$, $k=1,...,n$, then we have
\begin{equation}
\left\{ 
\begin{array}{l}
\int\limits_{\mathcal{O} }[\partial _{t}\mathbf{v}^{n}(P_{n}{\boldsymbol{\varphi }
})-\mathbf{v}^{n}\otimes \mathbf{v}^{n}:\nabla (P_{n}{\boldsymbol{\varphi }}
)+\mathbb{S}^{n}:\nabla (P_{n}{\boldsymbol{\varphi }})]\,d{x}\\ 
\qquad \qquad \qquad +\int_{\Gamma }(\nabla g(\mathbf{v}^{n})\cdot (P_{n}{
\boldsymbol{\varphi }}))\,d{\Gamma }\,dt=0, \\ 
\mathbf{v}^{n}(0)=\mathbf{v}_{0}^{n},
\end{array}
\right.  
\label{vnn}
\end{equation}%
since $\{\mathbf{e}_{j}\}_{j=1}^{\infty }$ \ is the orthogonal basis\ for
the space $W.$ By the H\"{o}lder inequality, the Sobolev continuous
embedding $W^{2,2}(\mathcal{O} )\subset C(\overline{\mathcal{O} })$ and the condition \eqref{condition for g}$_2$ for the boundary function $g$, we derive 
\begin{eqnarray*}
|(\partial _{t}\mathbf{v}^{n},{\boldsymbol{\varphi }})| &\leqslant &\left(
\Vert \mathbf{v}^{n}\Vert _{L^{2}(\mathcal{O} )}^{2}+\Vert \mathbb{S}^{n}\Vert
_{L_{1}(\mathcal{O} )}\right) \Vert \nabla (P_{n}{\boldsymbol{\varphi }})\Vert
_{C(\overline{\mathcal{O} })} \\
&&+c\Vert\mathbf{v}^{n}\Vert _{L_{1}(\Gamma )}\Vert P_{n}\boldsymbol{\varphi }\Vert
_{C(\overline{\mathcal{O} })} \\
&\leqslant &C\left(\Vert \mathbf{v}^{n}\Vert _{L^{2}(\mathcal{O} )}^{2}+\Vert 
\mathbb{S}^{n}\Vert_{L^{p^{\prime}}(\mathcal{O} )}+\Vert\mathbf{v}^{n}\Vert_{L^p(\Gamma)}\right)\Vert {P_{n}\boldsymbol{\varphi}}\Vert_{W}.
\end{eqnarray*}
The a priori estimates \eqref{omegadelta}-\eqref{Sdelta} and the continuous embedding
$V_p\subset H^1(\mathcal{O})$ together with the continuity of the trace map $V_p\equiv W^{1,p}(\mathcal{O})\to L^p(\Gamma)$ imply 
\begin{equation*}
\int_{0}^{T}|(\mathbf{v}^{n},\partial _{t}{\boldsymbol{\varphi }}%
)|dt\leqslant C\Vert \boldsymbol{\varphi }\Vert _{L^{p}(0,T;W)},
\end{equation*}
which gives \eqref{timedelta}. $\hfill\;\blacksquare$

\bigskip

\section{Limit transition}
\label{Limit}

The estimates \eqref{omegadelta}-\eqref{Sdelta}, together with the continuity of the trace map $V_p\to L^p(\Gamma)$ and the condition \eqref{condition for g}$_2$ for the boundary function $g$, imply the existence of subsequences satisfying the conditions \eqref{conv introducao}.

\smallskip

Let us remember the Aubin-Lions-Simon compactness result (see \cite{sim,tem}).

\begin{lemma}
\label{ALS}
Let $X_{0}$, $X$ and $X_{1}$ be three Banach spaces with $%
X_{0}\subseteq X\subseteq X_{1}$. Suppose that $X_{0}$, $X_{1}$ are
reflexive, $X_{0}$ is compactly embedded in $X$ and that $X$ is continuously
embedded in $X_{1}$. Let 
\begin{equation*}
\mathcal{Z}=\left\{ v\in L^{2}(0,T;X_{0}),\qquad \partial _{t}v\in
L^{p^{\prime }}(0,T;X_{1})\right\} .
\end{equation*}%
Then the embedding of $\mathcal{Z}$ into $L^{2}(0,T;X)$ is compact.
\end{lemma}

Hence the estimates \eqref{omegadelta}, \eqref{timedelta}, the compact
embedding $W^{1,p}(\mathcal{O} )\subset L^{2}(\mathcal{O} )$, the continuous
embedding $L^{2}(\mathcal{O})\subset W^{\ast }$ and Lemma \ref{ALS} give
the strong convergence $\mathbf{v}^{n}\rightarrow \mathbf{v}$ in $L^{2}(\mathcal{O}_{T})$, as stated previously in \eqref{convergence introducao}.

Applying the above convergences (more precisely, \eqref{conv introducao}, \eqref{convergence introducao} and \eqref{c02}) in \eqref{vnn}, we deduce that the limit functions $\mathbf{v}$, $\mathbb{S}$, $\mathbf{g}$ \ fulfil the integral equality
\begin{equation}
\int\limits_{\mathcal{O}_{T}}\left[\mathbf{v}\partial_{t}{\boldsymbol{\varphi}}+\left(\mathbf{v\otimes v}-\mathbb{S}\right):\nabla {\boldsymbol{\varphi}}\right]\,d{x}dt+\int\limits_{\mathcal{O}}\mathbf{v}_{0}{\boldsymbol{\varphi}}(x,0)\,d{x}=\int_{\Gamma_{T}}\mathbf{g}\cdot {\boldsymbol{\varphi}}\,d{\Gamma}\,dt  
\label{weak}
\end{equation}
for any function $\boldsymbol{\varphi }\in C^{1}([0,T];W)$, \ $\boldsymbol{\varphi }(\cdot,T)=0$.

\bigskip

To demonstrate the relation \eqref{SR01g} and that $\mathbf{g}=\nabla g(\mathbf{v})$  we need the condition $p\ge 2.2$. Indeed, with this condition we can prove the following result.

\begin{lemma}
\label{1 copy(1)} 
Under the conditions (hypothesis) stated in Theorem \ref{theorem 1} - in particular we stress the condition $p\ge 2.2$ -, \ the function 
\begin{equation*}
\mathbf{v}\in L^{\infty }(0,T;H)\cap L^{p}(0,T;V_{p})
\end{equation*}
fulfilling the integral equality \eqref{weak} satisfies the estimate
\begin{equation}
\Vert \partial_{t}\mathbf{v}\Vert_{L^{p^{\prime }}(0,T;\,\,V_{p}^{\ast
})}\leqslant C,  
\label{timedelta0}
\end{equation}
for some positive constant $C$ (that may depend on $p,$ $q,$ $\mathbf{v}_{0},$ $\mu _{i},$ ${\nu} $\ and the boundary function $g$).
\end{lemma}

\textbf{Proof}.\ \ Let ${\boldsymbol{\varphi }}\in C_{0}^{1}(0,T;W)$ be an
arbitrary function in \eqref{weak}. In what follows we consider $2\leqslant p<3$. The case $p\geqslant 3$ is a more simple and can be done by the same manner. Applying the H\"{o}lder inequality to \eqref{weak}, we derive 

\begin{eqnarray*}
\int\limits_{0}^{T}|(\mathbf{v},\partial _{t}{\boldsymbol{\varphi }})dt|
&\leqslant &\int\limits_{0}^{T}\{\Vert \mathbf{v}\Vert _{L^{\frac{3p}{4p-6}%
}(\mathcal{O} )}\Vert \nabla \mathbf{v}\Vert _{L^{p}(\mathcal{O} )}\Vert {\boldsymbol{%
\varphi }}\Vert _{L^{\frac{3p}{3-p}}(\mathcal{O})} \\
&&+\Vert \mathbb{S}\Vert_{L^{p^{\prime}}(\mathcal{O})}\Vert\nabla{\boldsymbol{\varphi }}\Vert_{L^{p}(\mathcal{O})} + \Vert \mathbf{g}\Vert_{L^{2}(\Gamma
)}\Vert {\boldsymbol{\varphi}}\Vert_{L^{2}(\Gamma)}\}dt \\
&\equiv&\int\limits_{0}^{T}\left\{ I_{1}+I_{2}+I_{3}\right\} dt.
\end{eqnarray*}

Using that $\mathbb{S}\in L^{p'}(\mathcal{O})$, $\mathbf{g}\in L^p(\Gamma_T)$ and the continuity of the trace map $W^{1,p}(\mathcal{O})\to L^{2}(\Gamma)$, we easily obtain
\begin{equation}
\int_{0}^{T}(I_{2}(t)+I_{3}(t))\,dt\leqslant C\Vert \boldsymbol{\varphi}\Vert _{L^{p}(0,T;V_{p})}.  
\label{i12}
\end{equation}

By Lemma \ref{equivalents}, the embedding \eqref{embed}\ and the interpolation inequality \eqref{interp} for $r=\frac{3p}{4p-6},$ $s_{1}=2,$ $s_{2}=\frac{3p}{3-p},$\ $\theta _{1}=1-\theta _{2}$, $\ \theta _{2}=\frac{12-5p}{5p-6}$, we have 
\begin{align}
\int_{0}^{T}I_{1}(t)\,dt\leqslant & \int_{0}^{T}\Vert \mathbf{v}\Vert _{L^{%
\frac{3p}{4p-6}}(\mathcal{O} )}\Vert \nabla \mathbf{v}\Vert _{L^{p}(\mathcal{O}
)}\Vert {\boldsymbol{\varphi }}\Vert _{L^{\frac{3p}{3-p}}(\mathcal{O} )}\,dt 
\notag \\
\leqslant & C\int_{0}^{T}\{\Vert \mathbf{v}\Vert _{L^{2}(\mathcal{O} )}^{\theta
_{1}}\Vert \nabla \mathbf{v}\Vert _{L^{p}(\mathcal{O} )}^{\theta _{2}}+\Vert 
\mathbf{v}\Vert _{L^{2}(\mathcal{O} )}\}\Vert \nabla \mathbf{v}\Vert
_{L^{p}(\mathcal{O} )}\Vert \boldsymbol{\varphi }\Vert _{V_{p}}\,dt  \notag \\
\leqslant & C\{\Vert \mathbf{v}\Vert _{L^{\infty }(0,T;L^{2}(\mathcal{O}
))}^{\theta _{1}}\left( \int_{0}^{T}\Vert \nabla \mathbf{v}\Vert
_{L^{p}(\mathcal{O} )}^{(1+\theta _{2})p^{\prime }}\,dt\right) ^{\frac{1}{%
p^{\prime }}}  \notag \\
& +\Vert \mathbf{v}\Vert _{L^{\infty }(0,T;L^{2}(\mathcal{O} ))}\left(
\int_{0}^{T}\Vert \nabla \mathbf{v}\Vert _{L^{p}(\mathcal{O} )}^{p^{\prime
}}\,dt\right) ^{\frac{1}{p^{\prime }}}\}\Vert \boldsymbol{\varphi }\Vert
_{L^{p}(0,T;V_{p})}  \label{vvn}
\end{align}%
by the H\"{o}lder inequality. We have 
\begin{equation*}
p^{\prime }\leqslant p,\quad (1+\theta _{2})p^{\prime }\leqslant p\quad 
\text{if}\quad p\geqslant 2.2,
\end{equation*}%
therefore combining \eqref{i12} and \eqref{vvn} with the help of the a
priori estimates \eqref{omegadelta}-\eqref{Sdelta}, we deduce 
\begin{equation*}
|\int_{0}^{T}(\mathbf{v},\partial _{t}{\boldsymbol{\varphi }})dt|\leqslant
C\Vert \boldsymbol{\varphi }\Vert _{L^{p}(0,T;V_{p})},
\end{equation*}%
which implies \eqref{timedelta0}. \ $\blacksquare$

\bigskip

To demonstrate the relation \eqref{SR01g} we use the approach of the theory of variational inequalities \cite{DL,Evans}. This was done in details in the Introduction and ends the proof of Theorem \ref{theorem 1}. 

\medskip
\textbf{Acknowledgment}  {\footnotesize {\ The authors would like to thank Prof. Anderson L. A. de Araujo (UFV, Brazil) for
useful discussions under the article subject.}

{\ The work of N.V. Chemetov was
supported by FAPESP (Funda\c{c}\~{a}o de Amparo \`{a} Pesquisa do Estado de S%
\~{a}o Paulo), project 2021/03758-8, "Mathematical problems in fluid
dynamics" and the Visiting Professor Projects
2023/11369-7, ""Fluid-rigid body" interaction problem", 2023/11593-4, ""Non-newtonian fluid - Rigid body" interaction problem". }}


\begin{thebibliography}{99}
\bibitem{AKM} \textsc{Antontsev S.N., Kazhikhov A.V., Monakhov V.N.,} \emph{%
Boundary value problems in mechanics of nonhomogeneous fluids.} Elsevier,
Amsterdam, 1990.

\bibitem{BC} \textsc{Bauschke H.H., Combettes L.}, \emph{Convex Analysis and
Monotone Operator Theory in Hilbert Spaces}. Springer (2011).



\bibitem{Boyer} \textsc{Boyer F., Fabrie P.,} \emph{Mathematical Tools for
the Study of the Incompressible Navier-Stokes Equations and Related Models.}
Springer Science+Business Media New York, 2013.

\bibitem{Brezis} \textsc{Brezis, H.,} (2011). \emph{Functional analysis, Sobolev Spaces and Partial Differential Equations}. New York: Springer-Verlag, 2011.



\bibitem{BR}\textsc{Bridges S., Robinson L.,} \emph{A Practical Handbook for Drilling Fluids Processing.} Elsevier (Imprint: Gulf Professional Publishing) 2020. DOI https://doi.org/10.1016/C2019-0-00458-X.

\bibitem{C} \textsc{Chemetov N.V., Antontsev S.N.}, \emph{Euler equations
with non-homogeneous Navier slip boundary condition}, Physica D: Nonlinear
Phenomena 237 (2008), pp. 92-105.

\bibitem{CC3} \textsc{Chemetov  N.V., Cipriano F.}, \emph{Boundary layer
problem: Navier-Stokes equations and Euler equations}, Nonlinear Analysis:
Real World Applications 14 (6) (2013), pp. 2091-2104.

\bibitem{CC5} \textsc{Chemetov N.V., Cipriano F.}, \emph{Inviscid limit for
Navier--Stokes equations in domains with permeable boundaries}, Applied
Math. Letters 33 (2014), pp. 6-11.

\bibitem{CC6} \textsc{Chemetov  N.V., Cipriano F., Gavrilyuk S.}, \emph{%
Shallow water model for the lake with friction and penetration},
Mathematical Methods in the Applied Sciences 33 (6) (2010), pp. 687-703.

\bibitem{Anderson-Nikolai} \textsc{Chemetov N.V., de Araujo A.L.A.}, \emph{
Well-posedness of the Cosserat-Bingham fluid equations}, NoDEA Nonlinear Differential Equations Appl. 29 (2022), no. 3, Paper No. 31, 24 pp.




\bibitem{Demengel} \textsc{Demengel, F., Demengel, G.,} \emph{Functional
spaces for the theory of elliptic partial differential equations.}
Translated from the 2007 French original by Reinie Ern\'e. Universitext.
Springer, London; EDP Sciences, Les Ulis, 2012.

\bibitem{DL} \textsc{Duvaut G., Lions J.L.,} \emph{Inequalities in Mechanics
and Physics.} Grundlehren der Mathematischez Wissenschaften, \textbf{212},
Springer-Verlag, Berlin, New York, 1976.

\bibitem{E} \textsc{Eringen A.C.,} \emph{Microcontinuum Field Theories.} \textbf{I}, \textbf{II}, Springer-Verlag, Berlin, New York, 1999.

\bibitem{Evans} \textsc{Evans L.C.,} \emph{Partial Differential Equations.} 
\textit{American Mathematical Society,} 1998.



\bibitem{G} \textsc{Gajewski F., Groger K., Zacharias K.} \emph{Nichtlineare
Operatorgleichungen und Operatordifferentialgleichungen.} Akademie-Verlag,
Berlin, 1974.

\bibitem{giusti} \textsc{Giusti E.} \emph{Direct Methods in the Calculus of Variations}. World Scientific Publishing Co. Pte. Ltd., 2003.






\bibitem{Lukas} \textsc{{\L}ukaszewicz G.}, \emph{Microploar Fluids. Theory and Applications}, Birkh\"auser Boston Inc., Boston, MA (1999).

\bibitem{Malek-Necas-Rokyta-Ruzicka} \textsc{Málek, J., Ne\v{c}as, J., Rokyta, M., R\r{u}\v{z}i\v{c}ka M.}, {\em Weak and Measure-valued Solutions to Evolutionary PDEs}. Publisher Chapman \& Hall, London (1996).





\bibitem{R} \textsc{Rockafellar R. T.}, \emph{Convex Analysis}, Princeton
University Press (1970).

\bibitem{smr} \textsc{R\r{u}\v{z}i\v{c}ka M., Shelukhin V.V., Santos M.M.,} 
\emph{Steady flows of Cosserat--Bingham fluids.} \ Math. Methods in the
Applied Sciences, \textbf{40}, Issue 7 (2017), 2746-2761.

\bibitem{ShelChem} \textsc{Shelukhin V.V., Chemetov N.V.,} \emph{Global
Solvability of the One-Dimensional Cosserat Bingham Fluid Equations.} J.
Math. Fluid Mech., \textbf{17}, No. 3 (2015), 495-511.


\bibitem{ShelRuz2013} \textsc{Shelukhin V.V., R\r{u}\v{z}i\v{c}ka M.,} \emph{%
On Cosserat--Bingham Fluids.} Z. Angew. Math. Mech., \textbf{93}, No. 1
(2013), 57--72.

\bibitem{T} \textsc{Teisseyre R., Teisseyre-Jerenska M.,} \emph{Asymmetric continuum: Extreme processes in solids and fluids.} Springer, Heldelberg, New York, London, 2014.

\bibitem{sim} \textsc{Simon J.,} \emph{Compact sets in the space} $%
L^{p}(0,T;B).$ Ann. Mat. Pura Appl., \textbf{146} (1987) 65-96.

\bibitem{tem} \textsc{Temam R.,} \emph{Navier-Stokes equations, Theory and
Numerical analysis.} AMS\ Chelsea Publishing, Providence, Rhode Island, 2001.

\bibitem{tem2} \textsc{Temam R.,} \emph{Mathematical Problems in Plasticity.} Gauthier-Villars, 1985.

\end{thebibliography}
\end{document}